\documentclass[12pt]{amsart} 

\usepackage{color} 
\usepackage{amsfonts} 
\usepackage{amsthm}
\usepackage{amsmath}
\usepackage{amssymb}
\usepackage{latexsym}
\usepackage[all]{xy}

\newcommand{\bF}{{\bf F}}
\newcommand{\bG}{{\bf G}}
\newcommand{\M}{\mathbb{M}}

\newcommand{\A}{\mathbb{A}}
\newcommand{\B}{\mathbb{B}}
\newcommand{\C}{\mathcal{C}}
\newcommand{\D}{\mathcal{D}}

\newcommand{\Z}{\mathbb{Z}}

\newcommand{\nF}{F}
\newcommand{\nG}{G}

\newcommand{\Hom}{{\rm Hom}}
\newcommand{\oHom}{{\rm \overline Hom}}
\newcommand{\End}{{\rm End}}

\newcommand{\Obj}{{\rm Obj}}
\newcommand{\Mor}{{\rm Mor}}

\newcommand{\vareps}{\varepsilon}
\newcommand{\rro}{\varrho}

\newcommand{\id}{I}

\newcommand{\cC}{{\mathcal C}}
\newcommand{\cD}{{\mathcal D}}

\newcommand{\kB}{{[B,-]}}

\newcommand{\kC}{{[\cC,-]}}

\newcommand{\oB}{{\overline B}} 
\newcommand{\uB}{{\underline B}}

\newcommand{\oH}{{\overline H}} 
\newcommand{\uH}{{\underline H}}
\newcommand{\oT}{{\overline T}}

\newcommand{\wT}{{\widehat T}}

\newcommand{\ot}{\otimes}
\newcommand{\lra}{\longrightarrow}
\newcommand{\LRa}{\Leftrightarrow}
\newcommand{\Ra}{\Rightarrow}

\newcommand{\roM}{\varrho^M}

\newcommand{\roN}{\varrho^N}
\newcommand{\Nro}{{^N\negthickspace\varrho}}

\renewcommand{\AA}{A}

\newcounter{zlist}
\newenvironment{zlist}{\begin{list}{{\rm(\arabic{zlist})}}{
\usecounter{zlist}\leftmargin2.5em\labelwidth2em\labelsep0.5em
\topsep0.6ex\itemsep0.3ex plus0.2ex minus0.3ex
\parsep0.3ex plus0.2ex minus0.1ex}}{\end{list}}
 
\newcounter{blist}
\newenvironment{blist}{\begin{list}{{\rm(\alph{blist})}}{
\usecounter{blist}\leftmargin2.5em\labelwidth2em\labelsep0.5em
\topsep0.6ex \itemsep0.3ex plus0.2ex minus0.3ex
\parsep0.3ex plus0.2ex minus0.1ex}}{\end{list}}

\newcounter{rlist}
\newenvironment{rlist}{\begin{list}{{\rm(\roman{rlist})}}{
\usecounter{rlist}\leftmargin2.5em\labelwidth2em\labelsep0.5em
\topsep0.6ex\itemsep0.3ex plus0.2ex minus0.3ex
\parsep0.3ex plus0.2ex minus0.1ex}}{\end{list}}

 \def\rhoma#1#2{{{\rm Hom}_A\left(#1,#2\right)}}
\def\ota{\otimes_A}
 \def\DC{{\Delta}} 
 
\def\rhom#1#2#3{{{\rm Hom}\sb{#1}(#2,#3)}}
\def\sw#1{{\sb{(#1)}}} 

\def\otc{{\otimes_{[\cC,-]}}}
\def\Rhom#1#2#3{{{\rm Hom}\sp{#1}\left(#2,#3\right)}} 

\def\Rhom#1#2#3{{{\rm Hom}\sp{#1}\left(#2,#3\right)}}

\def\coten#1{\ot^{#1}}

\def\can{{\rm \textsf{can}}} 
\def\tw{{\rm \textsf{tw}}} 

\swapnumbers
\newtheorem{theorem}{Theorem}[section]
\newtheorem{lemma}[theorem]{Lemma}

\newtheorem{corollary}[theorem]{Corollary}

\newtheorem{thm}[theorem]{} 

\def\sw#1{{\sb{\underline{#1}}}}

\def\Label{\label}

\def\arg#1{#1}             
\def\narg#1{#1}            
\def\CGbcm{$[\C,-]$-Galois bicomodule}
\def\DGbcm{$\D$-Galois bicomodule}
\def\CGcm{$[\C,-]$-Galois comodule}
\def\CGlcm{$[\C,-]$-Galois left comodule}
\def\DGlcm{$\D$-Galois left comodule}
\def\DGrcm{$\D$-Galois right comodule}
\def\counit{\varepsilon}
\def\cop{\Delta}
\def\product{\mu}
\def\unit{\iota}
\def\adjcu{\varepsilon}
\def\iso{\cong}
\def\equi{\simeq}

\title{Monads and comonads in module categories}
 \author[G. B\"ohm]{Gabriella B\"ohm}
 \address{Research Institute for Particle and Nuclear Physics, Budapest, 
 \newline\indent H-1525
 Budapest 114, P.O.B.\ 49, Hungary}
  \email{G.Bohm@rmki.kfki.hu}
 \author[T. Brzezi\'nski]{Tomasz Brzezi\'nski}
  \address{ Department of Mathematics, Swansea University, 
  Singleton Park, \newline\indent  Swansea SA2 8PP, U.K.} 
  \email{T.Brzezinski@swansea.ac.uk} 
\author[R. Wisbauer]{Robert Wisbauer}
\address{Department of Mathematics, Heinrich-Heine University, \newline\indent
D-40225
D\"usseldorf, Germany}
 \email{wisbauer@math.uni-duesseldorf.de}   
   \date{April 2008}
  
\begin{document}

 \begin{abstract}
Let $A$ be a   ring and $\M_A$ the category of $A$-modules. 
It is well known in module theory that  
for any $\AA $-bimodule $B$, {\em $B$ is an $A$-ring} if and only if the
functor $-\otimes_A B: \M_A\to \M_A$ is a {\em monad} (or {\em triple}). 
 Similarly, an $\AA $-bimodule $\C$ is an {\em $A$-coring}  
provided the functor
$-\otimes_A\C:\M_A\to \M_A$ is a {\em comonad} (or {\em cotriple}). 
The related categories of {\em modules} (or {\em algebras}) of $-\otimes_A B$
and {\em comodules} (or {\em coalgebras}) of $-\otimes_A\C$ are well studied 
in the literature. On the other hand, the right adjoint endofunctors
$\Hom_A(B,-)$ and $\Hom_A(\C,-)$ are a comonad and a monad, respectively, but  
the corresponding (co)module  categories did not find much attention so far. 
The category of $\Hom_A(B,-)$-comodules is isomorphic to the 
category of $B$-modules, while the category of $\Hom_A(\C,-)$-modules (called 
{\em $\C$-contramodules} by Eilenberg and Moore) need not be equivalent 
to the category of $\C$-comodules. 

The purpose of this paper is to investigate these categories and their 
relationships based on some observations of the categorical background.  
This leads to a deeper understanding and characterisations
of algebraic structures such as corings, bialgebras and Hopf algebras. 
For example, it turns out that the categories of $\C$-comodules and 
$\Hom_A(\C,-)$-modules are equivalent
provided $\C$ is a coseparable coring. Furthermore, a bialgebra
 $H$ over a commutative ring $R$ is a Hopf algebra 
if and only if $\Hom_R(H-)$ is a 
Hopf bimonad on $\M_R$ and in this case the categories of $H$-Hopf modules
and mixed $\Hom_R(H,-)$-bimodules are both equivalent to $\M_R$.  
\end{abstract}
 \maketitle
 
\tableofcontents 

\section{Introduction}

The purpose of this paper is to present a categorical framework for studying 
problems in the theories of rings and modules, corings and comodules, bialgebras
and (mixed) bimodules and Hopf algebras and Hopf modules.
The usefulness of this framework is illustrated by analysing the structure of
the  
category of {\em contramodules} and the bearing of this
structure on the properties of corings and 
bialgebras.

It is well-known that for a right module $V$ over an $R$-algebra 
$A$, the dual $R$-module $V^* = \rhom R V R$ is a left module over $A$. 
It is equally well-known that for a right comodule $V$ of an
$R$-coalgebra $C$, in general $V^*$ is not a $C$-comodule (left or right). 
It has already been 
realised in \cite[Chapter IV.5]{EilMoo:fou} that to a coalgebra $C$ 
two (different) representation categories can be associated: the familiar 
category of $C$-comodules and the category of {\em $C$-contramodules}
 introduced therein. If $V$ is a $C$-comodule, then $V^*$ 
 is a $C$-contramodule.

While comodules of coalgebras (and corings) have been 
intensively studied,  contramodules seem to have been rather neglected.
Yet the category of contramodules is  as fundamental as that of comodules, 
and  both categories are complementary to each other. To substantiate
this claim, one needs to resort to the categorical point of view on corings.
An $A$-coring can 
be defined as an $\AA $-bimodule $\cC$ such that the tensor endofunctor 
$-\ot_A\cC$ on the category of right $A$-modules $\M_A$ is a comonad or a
cotriple. Right $\cC$-comodules are the same as {\em comodules} (or {\em
  coalgebras} in category theory terminology)
of the comonad $-\ot_A\cC$. On the other hand, the tensor functor 
$-\ot_A\cC$ has a right adjoint, the Hom-functor $\rhom A \cC -$. 
By purely categorical arguments (see Eilenberg 
and Moore \cite[Proposition 3.1]{EM}), the functor $-\ot_A\cC$ is a 
{\em comonad} if and only if its {\em right} adjoint $\rhom A \cC -$ 
is a {\em monad}. Thus, $\cC$ is an $A$-coring  
if and only if $\rhom A \cC -$ is a monad on $\M_A$; right 
$\cC$-contramodules are simply {\em modules} (or {\em algebras}
in category theory terminology) of this monad. 
This categorical interpretation explains the way in which contramodules 
complement comodules. For example, since $\cC$-comodules are 
comodules of a comonad on an abelian category, their category has 
cokernels but not necessarily kernels. On the other hand, 
since $\cC$-contramodules are modules of a monad, their category 
has kernels but not necessarily cokernels. Thus one category 
provides the structure which the other one misses. 

Again purely categorical considerations (see \cite{EM}) explain that, while 
there are two categories of representations of a coring, there is only one
category of representations of a ring -- the familiar category of modules.
More precisely, a ring morphism $A\to B$ can be equivalently
described as the monad structure of the tensor functor $-\ot_A B$ on 
$\M_A$ associated to an $\AA $-bimodule $B$. With this interpretation,
right $B$-modules are simply modules of the monad $-\ot_A B$. 
The right adjoint functor $\rhom A B -$ is a comonad on 
$\M_A$ and the category of comodules of $\rhom A B -$ is {\em isomorphic} 
to the category of modules of the monad $-\ot_A B$. 
Consequently, there is only one type of representation categories for rings -- 
the category of right (or left) modules over a ring. Since modules of a ring
are thus both algebras and coalgebras of respective monads and comonads, the
category of modules inherits both kernels and cokernels from the category of
abelian groups. 

The above comments illustrate how the  categorical point of view can give  
significant insight into algebraic structures. There are many constructions
developed in category theory that are directly applicable to ring
theoretic situations  but they seem not to be sufficiently
explored. Contramodules of a coring are a good example of this. On one hand,
from the category point of view, they are as natural as comodules, on the
other hand, their structure was not analysed properly until very recently,
when their important role in  semi-infinite homology
was outlined by Positselski \cite{Pos}.  
The main motivation of our paper is a study of contramodules of corings. This
aim is achieved by placing it in a broader context: we revisit category
theory, more specifically the theory of adjoint comonad-monad pairs, in the
context of rings and modules. 
 
We begin by summarising the categorical framework, and then apply it first to
rings in module categories, next to corings. In the latter case, we concentrate
on properties of the less-known category of 
contramodules, and derive consequences of the 
categorical formulation in this context. We analyse functors between 
categories of comodules and contramodules, and study equivalences between
such categories involving a Galois theory for bicomodules. 
We also derive 
the characterisation of entwining structures as liftings of Hom-functors to
module and contramodule categories. 

Finally, we study contramodules of corings 
associated to bialgebras and provide new extensions of the Fundamental
Theorem of Hopf algebras (see \ref{Hopf.char}). First we observe that an 
$R$-module $B$ over a commutative ring $R$ is a bialgebra if and only if 
$\Hom_R(B,-)$ is a {\em bimonad}, that is, a monad and a comonad on $\M_R$
satisfying some compatibility conditions (see \ref{bimonad.B}). The
Fundamental Theorem says that a bialgebra $B$ is a Hopf algebra if and only if
$-\ot_RB$ induces an equivalence
between $\M_R$ and the category $\M^B_B$ of Hopf modules.
This can also be formulated as $B$ being a Galois comodule of associated
corings. 
Here we add that a Hopf algebra $B$ is characterised by a bimonad $\Hom_R(B,-)$ 
inducing an equivalence between $\M_R$ and the category $\M^\kB_\kB$ of certain 
$\Hom_R(B,-)$-bimodules  ({\em Hopf contramodules}).
 Again this can be seen as $B$ being a {\em Galois comodule} 
with respect to the Hom-functors of the associated corings. 

\bigskip

 \section{Categorical framework}

Our main concern is to apply abstract categorical notions to 
special situations in module categories. We begin by recalling
some basic definitions and properties (e.g.\ from \cite{EM}) 
to fix notation, and then develop a categorical framework 
which is later applied to categories of (co)modules. 

Throughout, the composition of functors is denoted by juxtaposition, and the 
usual composition symbol $\circ$ is reserved for natural transformations and
morphisms. Given functors $F$, $G$ and a natural transformation $\varphi$, 
$F\varphi G$ denotes the natural transformation, which, evaluated at an object
$X$ gives a morphism obtained by applying $F$ to a morphism provided by the
natural transformation $\varphi$ evaluated at the object $G\arg X$.

By $\A\equi \B$ we denote equivalences between categories and 
$\A\iso \B$ is written for their isomorphisms. The symbol
$\iso$ is also used to denote isomorphisms between objects in any 
category, in particular isomorphisms of modules and (natural) isomorphisms of
functors.

 \begin{thm}{\bf Adjoint functors.} \em
 A pair $(L,R)$ of functors   
  $L:\A\to \B$ and $R:\B\to \A$ between categories $\A$, $\B$   
 is called an {\em adjoint pair} if there is a natural isomorphism 
   $$ \Mor_\B (L(-) ,-) \stackrel{\iso}\lra \Mor_\A (-,R(-)). $$
This can be described by 
  natural transformations  {\em unit} $\eta:\id_\A \to RL$
and {\em counit}  $\adjcu: LR\to \id_\B$
 satisfying the 
    {\em triangular identities}
$\adjcu L\circ L\eta = \id_L$ and $R\adjcu \circ \eta R = \id_R$.
 
Recall the properties    
 of an adjoint pair $(L,R)$: 
{\em  \begin{zlist}
  \item $R$ is full and faithful if and only if $\adjcu:LR\to \id_\B$
    is an isomorphism.
  \item $L$ is full and faithful if and only if $\eta: \id_\A\to RL $ 
     is an isomorphism.
  \item $L$ is an equivalence if and only if $\adjcu$ and $\eta$ are
    isomorphisms.  
\end{zlist} } 
\end{thm}  
  
\begin{thm}\Label{nat.adj}{\bf Natural transformations for adjoints.} \em 
For two adjunctions $(L,R)$ and $({\widetilde L},{\widetilde R})$ between $\A$
and $\B$, with
respective units $\eta$, ${\widetilde \eta}$ and counits $\adjcu$, 
${\widetilde \adjcu}$, there is an isomorphism between the natural
transformations (cf. \cite{KeSt}, \cite{MW}) 
$$
\mathrm{Nat}(L,{\widetilde L}) \to \mathrm{Nat}({\widetilde R}, R),\quad 
f \mapsto {\bar f}:= R {\widetilde \adjcu} \circ R f {\widetilde R}
\circ \eta {\widetilde R},
$$
with the inverse map
$$
\mathrm{Nat}({\widetilde R}, R) \to \mathrm{Nat}(L,{\widetilde L}),\qquad 
{\bar f} \mapsto f:=\adjcu {\widetilde L}\circ L{\bar f} {\widetilde
  L}\circ L {\widetilde \eta}. 
$$
We say that $f$ and ${\bar f}$ are {\em mates under the adjunctions
$(L,R)$ and $({\widetilde L},{\widetilde R})$.} For natural transformations 
$f:L_1 \to L_2$ and $g: L_2 \to L_3$ between left adjoint functors, naturality
and the triangle identities imply 
$\overline{g \circ f} = {\bar f} \circ {\bar g}$. In particular, $f$ is a
natural isomorphism if and only if its mate ${\bar f}$ is a natural
isomorphism. Moreover, if for an adjunction $(L,R)$, the composites $L L_1$
(and hence $L L_2$) are meaningful, then $\overline {Lf}={\bar f}
R$. Similarly, if the composites $L_1 L$ (and thus $L_2 L$) are meaningful
then $\overline{fL}=R {\bar f}$.
\end{thm}

 \begin{thm}\Label{monad}{\bf Monads on $\A$.} \em 
A {\em monad on the category $\A$} is a triple 
   $\bF= (F,m,i)$, where 
$F:\A\to \A$ is a functor with natural transformations 
$m:F  F\to F$ and $i: \id_\A \to F$ satisfying 
associativity and unitality conditions. 
A {\em morphism of monads} $(F,m,i)\to (F',m',i')$ is a natural transformation
$\varphi:F \to F'$ such that $m'\circ \varphi F'\circ F \varphi = \varphi\circ
m$ and $\varphi \circ i = i'$.

 An {\em $\nF$-module} is  
 a pair consisting of  $A \in\Obj(\A)$  and a morphism $\varrho_{A}: F\arg
 A\to A$  such that the following  diagrams
$$\xymatrix{
 F  F \arg A \ar[r]^{\;m \narg A}\ar[d]_{F\varrho_A} & F\arg
 A\ar[d]^{\varrho_A} &    
A\ar[r]^-{i \narg A\quad}\ar[rd]_-{\id_A} & F\arg A\ar[d]^{\varrho_A}   \\
   F \arg A\ar[r]_{\varrho_A} &\; A\, , & & \; A \,  } $$  
are commutative.

Morphisms between $F$-modules $f:A\to A'$ are morphisms in $\A$ 
with  $ \rro_{A'} \circ F \arg f= f\circ \rro_A $
and the {\em Eilenberg-Moore category of $\nF$-modules} is denoted by $\A_\nF$. 

 For any object $A$ of $\A$, $F\arg A$ is an $\nF$-module 
and this yields the {\em free functor}
  $$\phi_\nF: \A \to \A_\nF, \; A \mapsto (F\arg A, m \narg A), $$ 
which is left adjoint to the 
{\em forgetful functor} $U_\nF: \A_\nF\to \A$  by the isomorphism 
$$\Mor_{\A_\nF}(\phi_\nF \arg A, B) \to \Mor_\A( A, U_\nF \arg B), \quad
f\mapsto f\circ i \narg A.$$ 
The full subcategory of $\A_\nF$ consisting of all  free 
$\nF$-modules (i.e.\ the full subcategory of $\A_\nF$ generated by the 
image of $\phi_\nF$) is called the {\em Kleisli category of $\bF$} 
and is denoted by $\tilde\A_\nF$.
\end{thm}

 \begin{thm}\Label{comonad}{\bf Comonads on $\A$.} \em 
A {\em comonad on $\A$} is a triple $\bG=(G,d,e)$,  where 
$G:\A \to \A $ is a functor with natural transformations
$d:G\to G  G$ and $e: G\to \id_\A $ 
satisfying coassociativity and counitality conditions. 
A {\em morphism of comonads} is a natural transformation that is compatible
with the coproduct and counit.
 
 A {\em $\nG$-comodule} is an object $A\in \A$ with a morphism  
$\varrho^{A}:A\to G\arg A$ compatible with $d$ and $e$.
Morphisms between $\nG$-comodules $g:A\to A'$ are morphisms in $\A$ 
with  $\rro^{A'}\circ g= G\arg g \circ \rro^A$
and the {\em Eilenberg-Moore category of $\nG$-comodules} is denoted by
$\A^\nG$. 

For any $A\in \A$, $G\arg A$ is a $\nG$-comodule yielding the {\em (co)free
  functor} 
$$
\phi^\nG: \A \to \A^\nG, \; A \mapsto (G\arg A, d\narg A) 
$$
which is right adjoint to the 
forgetful functor $U^\nG:\A^\nG\to \A$ by the isomorphism 
 $$\Mor_{\A^\nG}(B,\phi^\nG \arg A)\to \Mor_{\A}(U^\nG\arg B,A),\quad f\mapsto
e \narg A\circ f.$$ 
The full subcategory of $\A^\nG$ consisting of all  (co)free 
$\nG$-comodules (i.e.\ the full subcategory of $\A^\nG$ generated by the 
image of $\phi^\nG$) is called the {\em Kleisli category of $\bG$} 
and is denoted by $\tilde\A^\nG$.
\end{thm}
    
 \begin{thm}\Label{mon-adj} {\bf (Co)monads related to adjoints.}  
Let $L:\A\to \B$ and $R:\B\to \A$ be an adjoint pair of functors
with unit $\eta:\id_\A \to RL$ and counit $\adjcu:LR\to \id_\B$. Then
  $$\bF=(RL,\; 
     RLRL \stackrel{R\adjcu L}\lra RL,\;  \id_\A\stackrel{\eta}\lra RL)$$ 
is a monad on $\A$. Similarly a comonad on $\B$ is defined by
 $$\bG = (LR,\; 
     LR \stackrel{L\eta R}\lra LRLR, \;LR\stackrel{\adjcu}\lra \id_\B).$$
\end{thm}

As already observed by Eilenberg and Moore in \cite{EM},
the monad structure of an endofunctor induces a comonad structure
on any adjoint endofunctor. 
  
 \begin{thm}\Label{adj-mon} {\bf Adjoints of monads and comonads.}  
Let $L:\A\to \A$ and $R:\A\to \A$ be an adjoint pair of functors
with unit $\eta:\id_\A \to RL$ and counit $\adjcu:LR\to \id_\A$. 
\begin{zlist}
\item The following are equivalent:
 \begin{blist}
 \item $L$ is a monad;
 \item $R$ is a comonad.
 \end{blist}
In this case the Eilenberg-Moore categories $\A_L$ and $\A^R$ are isomorphic
to each other. 

\item The following are equivalent:
 \begin{blist}
 \item $L$ is a comonad;
 \item $R$ is a monad.
 \end{blist}
 In this case the Kleisli categories $\tilde\A^L$ and $\tilde\A_R$ are
 isomorphic to each other.   
\end{zlist}

\end{thm}
\begin{proof}  
(1) This is shown in \cite[Proposition 3.1]{EM} and here it follows from 
\ref{nat.adj} (see also \cite{MW}).
 \smallskip

The isomorphism of $\A_L$ and $\A^R$ is mentioned in \cite[p.\ 3935]{Street}. 
For convenience we explain the relevant functor leaving objects and
morphisms unchanged but turning $L$-module structure maps to $R$-comodule 
structure maps and vice versa.
An $L$-module $\rro_A:L\arg A \to A$ induces a morphism
$$
\xymatrix{
 A\ar[r]^{\eta\narg A\quad} & RL\arg A \ar[rr]^{R \arg {\rro_A}}& & R\arg A,
}
$$ 
making $A$ an $R$-comodule. Similarly, a comodule 
  $\rro^A:A \to R\arg A$ induces
$$
\xymatrix{
L\arg A \ar[rr]^{L\arg {\rro^A}\quad}&& LR\arg A 
\ar[r]^{\quad\vareps\narg A} & A,}
$$
defining an $L$-module structure on $A$. 
In a word, $\varrho_A$ and $\varrho^A$ are mates under the adjunctions $(L,R)$
and $(\id_\A,\id_\A)$. 
\medskip

An $L$-module morphism $f:A'\to A$ yields a commutative diagram
$$
\xymatrix{ 
L\arg {A'} \ar[r]^{L\arg f} \ar[d]_{\rro_{A'}} 
& L\arg A \ar[d]^{\rro_A}  \\
A'\ar[r]^f & A \ , }$$
from which we obtain the commutative diagram
$$
\xymatrix{ 
A' \ar[r]^f \ar[d]_{\eta\narg {A'}} & 
A \ar[d]^{\eta\narg A} \\ 
RL \arg {A'} \ar[r]^{RL\arg f} \ar[d]_{R\arg {\rro_{A'}}} 
           & RL\arg A \ar[d]^{R\arg {\rro_A}}  \\
    R\arg {A'}\ar[r]^{R\arg f} & R\arg A .}$$
Commutativity of the outer rectangle shows that $f$ is also 
an $R$-comodule morphism.
Similarly one can prove that $R$-comodule morphisms are also
$L$-module morphisms.
\medskip

(2) The equivalence of (a) and (b) is proved similarly to (1).

The isomorphism of the Kleisli categories $\tilde\A^L$ and $\tilde\A_R$
was observed in \cite{Klein} and is also mentioned in \cite[p. 3935]{Street}.
It is provided by the 
canonical isomorphisms for $A,A'\in \A$,
$$\begin{array}{rcl}
\Mor_{\A^L}(\phi^L\arg A,\phi^L\arg {A'})&\iso& \Mor_\A(L\arg A,A') \\
   &\iso & \Mor_\A(A,R\arg {A'}) \iso \Mor_{\A_R}(\phi_R\arg A,\phi_R\arg
    {A'}). 
\end{array}$$
 \end{proof}

\begin{thm}\Label{rel.proj}
{\bf Relative projectivity and injectivity.} \em
An object $A$ of a category ${\mathbb A}$ is said to be {\em
  projective relative to a functor $F:{\mathbb A}\to {\mathbb B}$} (or
 {\em $F$-projective} in short) if $\Hom_{\mathbb A}(A,f):\Hom_{\mathbb A}(A,X)
\to \Hom_{\mathbb A}(A,Y)$ is surjective for all those morphisms $f$ in
${\mathbb A}$, for which $F\arg f$ is a split epimorphism in ${\mathbb
  B}$. Dually, $A\in 
{\mathbb A}$ is said to be {\em injective relative to $F$} (or {\em
  $F$-injective}) if $\Hom_{\mathbb A}(f,A):\Hom_{\mathbb A}(Y,A) \to
\Hom_{\mathbb A}(X,A)$ is surjective for all such morphisms $f$ in ${\mathbb
  A}$, for which $F\arg f$ is a split monomorphism in ${\mathbb B}$.

For an adjunction $(L:{\mathbb A} \to {\mathbb B}, R:{\mathbb B} \to {\mathbb
  A})$, with unit $\eta$ and counit $\adjcu$, an object $A\in {\mathbb A}$
is $L$-injective if and only if $\eta \narg A$ is a split monomorphism in
${\mathbb A}$. Dually, $B\in {\mathbb B}$ is $R$-projective if and only if
$\adjcu \narg B$ is a split epimorphism in ${\mathbb B}$.
\end{thm}

 Recall (e.g.\ from \cite[Section 6.5]{Bor:han1} or \cite[A.1, p.\ 62]{Laz})
 that the {\em Cauchy completion}, also called {\em Karoubian closure}, of 
any category ${\mathbb A}$ is the smallest category ${\overline {\mathbb A}}$
that contains ${\mathbb A}$ as a 
subcategory and in which idempotent morphisms split (i.e.\ can be written as a
composite of an epimorphism and its section). The Karoubian closure is
unique up to equivalence and can be constructed as follows. Objects of
${\overline{\mathbb A}}$ are pairs $(A,a)$, where $A$ is an object in
${\mathbb A}$ and $a:A\to A$ is an idempotent morphism (i.e.\ $a\circ
a=a$). Morphisms $(A,a)\to (A',a')$ are morphisms $f:A\to A'$ in 
${\A}$, such that $a'\circ f = f\circ a$. 
If $F:{\A}\to {\B}$ is
an isomorphism, then an isomorphism ${\overline{\A}}\to {\overline {\B}}$ 
is given by the object map $(A,a) \mapsto (F\arg A,F\arg a)$ and the
morphism map $f\mapsto F\arg f$.

 \begin{thm}\Label{Karoubi}
{\bf Equivalence of Karoubian closures.}  
Let $\A$ be a category in which idempotent morphisms split. Then
 for any comonad $(L,d,e)$ on $\A$ with right adjoint monad $(R,m,i)$, 
there is an equivalence
$$
E:  \A^L_{inj} \to \A_{R}^{proj},
$$
where $\A^L_{inj}$ denotes the full subcategory of $\A^L$ whose objects
are injective relative to the forgetful functor $U^L:\A^L\to \A$ and 
$\A_{R}^{proj}$ denotes the full subcategory of $\A_{R}$ whose
objects are projective relative to the forgetful functor 
$U_R:\A_{R}\to\A$. 

Explicitly, for $(A,\varrho^A)\in \A^L_{inj}$, the object
$E(A,\varrho^A)$ is given by the equaliser of the parallel morphisms 
$R\arg {\varrho^A}$ and $\omega:= m \narg {L \arg A} \circ R \arg{\eta \narg A}:
R\arg A \to RL\arg A$, where $\eta$ is the unit of the adjunction $(L,R)$.
\end{thm}

\begin{proof}
By \ref{adj-mon}(2), the Kleisli categories ${\tilde \A}^L$ and ${\tilde
  \A}_{R}$ are isomorphic. As recalled above, this isomorphism  
extends to their Karoubian closures. The Karoubian closure of ${\tilde
  \A}_{R}$ is equivalent to the 
full subcategory of $U_R$-projective objects of $\A_{R}$ (see \cite{GuRi},
\cite[Theorem 2.5]{RoWo}). Dually, the Karoubian closure of ${\tilde \A}^L$ is
equivalent to the full subcategory of $U^L$-injective objects of $\A^L$. This
proves the equivalence $\A^L_{inj}\equi \A^{proj}_{R}$. 

The explicit form of the equivalence functor is obtained by computing the
composite of the isomorphism between the Karoubian closures of the Kleisli
categories with the equivalences in 
\cite[Theorem 2.5]{RoWo}. This (straightforward) computation yields 
 the equaliser $E(A,\varrho^A)\to R\arg A$ of the identity morphism $\id_{R\arg
A}$ and the idempotent morphism $R \arg {\nu^A} \circ \omega :R\arg A \to
 R\arg A$,  
where $\nu^A$ is a retraction of $\eta(A,\varrho^A)= \varrho^A$ in 
${\mathbb A}^L$. This equaliser exists 
by the assumption that idempotents split in ${\mathbb A}$. 
Since
$$
\omega \circ R \nu^A \circ \omega = R \varrho^A \circ R \nu^A \circ \omega 
\quad \mbox{and} 
\quad R\arg{\nu^ A}\circ  R\arg{\varrho^ A} = \id_{RA},
$$
$E(A,\varrho^A)\to R\arg A$ is also 
an equaliser of $\omega$ and $R\arg{\varrho^ A}$.
\end{proof}

\bigskip

Recall from  \cite{NdO}, \cite{Rafael} that a functor
 $F:\B\rightarrow\A$ is said to be {\em separable}
 if and only if
the transformation
$
{\Mor_{\B}}(-,-)\rightarrow {\Mor_{\A}
}({F}(-),{F}(-)),$ $f\mapsto {F}\arg {f},
$
 is a split natural monomorphism. Separable functors reflect split
 epimorphisms and split monomorphisms.

Questions related to 2.9(1) are also discussed in 
\cite[Proposition 6.3]{BruVir}.

  \begin{thm}\Label{sep.mon} {\bf Separable monads and comonads.} 
 Let $\A$ be a category.
\begin{zlist}
\item For a monad $(R,m,i)$ on ${\A}$, the
  following are equivalent:
\begin{blist}
\item $m$ has a natural section ${\widehat m}$ such that
$$
 Rm \circ {\widehat m}R = {\widehat m} \circ m = mR \circ R{\widehat m} ;
$$
\item the forgetful functor $U_R:{\mathbb A}_R \to {\mathbb A}$ is
  separable.
\end{blist}
\item For a comonad $(L,d,e)$ on ${\A}$, the
  following  are equivalent: 
\begin{blist}
\item $d$ has a natural retraction ${\widehat d}$ such that
$$
 {\widehat d} L  \circ  Ld  = d\circ {\widehat d} =  L {\widehat d}  \circ  dL ;
$$
\item the forgetful functor $U^L:{\mathbb A}^L \to {\mathbb A}$ is
  separable. 
\end{blist}
\end{zlist}
\end{thm}

\begin{proof}
(1) By Rafael's theorem \cite[Theorem 1.2]{Rafael}, $U_R$ is separable 
if and only if the counit $\adjcu_R$ of the adjunction $(\phi_R,U_R)$ 
(see \ref{monad}) is a split natural epimorphism. 

(1) (a)$\Rightarrow$(b).
 A section $\nu: \id_{{\mathbb A}_R} \to \phi_R U_R$ of
$\adjcu_R$ is given by a morphism  
$$\xymatrix{
\nu (X,\rro_X)\ :\ 
X \ar[r]^-{i\narg X} & R\arg X \ar[r]^-{{\widehat m} \narg X}
& RR \arg X \ar[r]^-{R\arg \rro_X} & R\arg X ,}
$$
for any $(X,\rro_X)$ in ${\A}_R$.
By naturality and the properties of ${\widehat m}$ required in (a),
 $\nu (X,\rro_X)$ is an $R$-module morphism, i.e.\ 
$ m \narg X  \circ  R\arg{\nu (X,\rro_X)} =\nu (X,\rro_X) \circ \rro_X$. 
Since ${\widehat m}$ is a section of $m$, $\nu(X,\rro_X)$ is a
section of $\adjcu_R(X,\rro_X)=\rro_X$. In order to see that, use also
associativity and unitality of the $R$-action $\rro_X$.  
The morphism $\nu$ is natural, i.e.\ for $f:(X,\rro_X) \to
(Y,\rro_Y)$ in ${\mathbb A}_R$, $ R\arg f \circ \nu(X,\rro_X) =
\nu(Y,\rro_Y)\circ f$. This 
follows by definition of an $R$-module morphism and naturality.

(b)$\Rightarrow$(a). A section $\nu:\id_{{\A}_R} \to\phi_R U_R$ of 
$\adjcu_R$ induces a section of $m=U_R \adjcu_R \phi_R$ by putting
${\widehat m}:= U_R \nu \phi_R$. It obeys the properties in (a) by
naturality.\smallskip

(2) The proof is symmetric to (1).
\end{proof}

\begin{thm}\Label{adj-sep} 
{\bf Separability of adjoints.}  
Let $L:\A\to \A$ and $R:\A\to \A$ be an adjoint pair of functors
with unit $\eta:\id_\A \to RL$ and counit $\adjcu:LR\to \id_\A$. 

If $(L,d,e)$ is a comonad with corresponding  monad $(R,m,i)$,
then there are pairs of adjoint (free and forgetful) functors
(see \ref{monad}, \ref{comonad}):
  $$\begin{array}{rll}
\xymatrix{\A \ar[r]^{\phi_R} & \A_R},& 
\xymatrix{ \A_R \ar[r]^{U_R}& \A },&
\mbox{ with unit } \eta_R \mbox{ and counit } \adjcu_R, \mbox{ and} \\[+1mm] 
\xymatrix{\A^L \ar[r]^{U^L}  & \A}, &
\xymatrix{\A \ar[r]^{\phi^L}  & \A^L }, & 
\mbox{ with unit }\eta^L \mbox{ and counit }\adjcu^L.
\end{array}$$ 
\begin{zlist}
\item 
$\phi^L$ is separable if and only if $\phi_R$ is separable.
\item
$U^L$ is separable if and only if $U_R$ is separable.
\end{zlist}
If the properties in part (2) hold, then any object of ${\mathbb A}^L$ is
  injective relative to $U^L$ and every object of ${\mathbb A}_R$ is 
  projective relative to $U_R$.
\end{thm}

\begin{proof} 
(1) By Rafael's theorem \cite[Theorem 1.2]{Rafael}, 
$\phi^L$ is separable if and only if $\adjcu^L=e$ is a split
natural epimorphism, while $\phi_R$ is separable if and only if $\eta_R=i$ is a
split natural monomorphism. By construction, $i$ and $e$ are mates under the
adjunction $(L,R)$ and the trivial adjunction $(\id_{\mathbb A},\id_{\mathbb
  A})$.   
That is, $e=\adjcu\circ  Li $ equivalently, $i= Re \circ
\eta$. Hence a natural transformation ${\widehat i}: R \to \id_{\mathbb A}$ is a
retraction of $i$ if and only if its mate 
${\widehat e}:={\widehat i}L\circ\eta$ under the adjunctions 
$(\id_{\mathbb A},\id_{\mathbb A})$ and $(L,R)$ is a section of $e$. 
\smallskip

(2) Since $d$ and $m$ are mates under the adjunctions $(L,R)$ and
$(LL,RR)$, a natural transformation ${\widehat m}$ satisfies the properties in
\ref{sep.mon}(1)(a) if and only if its mate ${\widehat d}$ satisfies the
properties in \ref{sep.mon}(2)(a). Thus the claim follows by \ref{sep.mon}.
\smallskip

It remains to prove the final claims. 
Following \ref{rel.proj},  an $L$-comodule $(A,\varrho^A)$ is 
$U^L$-injective if and only if $\eta^L(A,\varrho^A)=\varrho^A$ is a split
monomorphism in ${\mathbb A}^L$. Since $\varrho^A$ is split in
${\mathbb A}$ (by $e\narg A$) and $U^L$, being separable, reflects split
monomorphisms, any $(A,\varrho^A)\in \A^L$ is $U^L$-injective.
$U_R$-projectivity of every object of ${\mathbb A}_R$ is proven by a
symmetrical reasoning.
\end{proof}

\begin{thm} \label{lift}    {\bf Lifting of functors.} \em 
 Let $F:\A\to \A$, $G:\B\to \B$ and $T:\A\to \B$ be functors on 
the categories $\A$, $\B$. If $F$, $G$ are monads or comonads,   
we may consider the diagrams 
 $$\xymatrix{ 
        \A_F \ar@{-->}[r]^\oT \ar[d]_{U_F} & \B_G \ar[d]^{U_G} \\
      \A \ar[r]^T &\B , }\qquad
    \xymatrix{ 
        \A^F \ar@{-->}[r]^\wT \ar[d]_{U^F} & \B^G \ar[d]^{U^G} \\
      \A \ar[r]^T &\B ,}   
  $$  
where the $U$'s denote the forgetful functors.
If there exist $\oT$  or $\wT$ making the corresponding diagram
 commutative, they are called {\em liftings of $T$}.
Their existence as well as their properties depend on (the existence of)
natural transformations between combinations of the functors involved. These 
are called {\em distributive laws}.
 \end{thm} 

The following theorem is proved by 
Johnstone in \cite[Lemma 1 and Theorem 2]{John}
(see also \cite[Section 3.3]{TTT}, \cite[3.3]{WisAlg}). 

  \begin{thm}\Label{Apple}  {\bf Lifting for monads.}  
Let $\bF=(F,m,i)$ and  $\bG=(G,m',i')$ be monads on the 
categories $\A$ and $\B$, respectively, and let $T: \A \to \B $ be a functor.

\begin{zlist}
\item  The liftings $\oT: \A_\nF\to \B_\nG$ of $T$ 
are in bijective correspondence 
with the natural transformations $\lambda: GT\to TF$ inducing  
 commutative diagrams
$$
\xymatrix{
GGT\ar[rr] ^{m' T} \ar[d]_{G \lambda} & & 
GT \ar[d]^\lambda \\
GTF \ar[r]^{\lambda F} & TFF \ar[r]^{Tm} & TF, }   \quad
    \xymatrix{ T\ar[r]^{i' T} \ar[dr]_{Ti} & GT \ar[d]^\lambda \\
          & TF.} $$
Such a pair $(T,\lambda)$ is called a {\em  monad 
morphism} in \cite{LaSt}.

\item  If $T$ has a left adjoint and $\A$ has coequalisers of reflexive 
  $TU_\nF$-contractible coequaliser pairs,
 then any lifting $\oT$ has a left adjoint.	  
\end{zlist}   
\end{thm}

For endofunctors the preceding diagrams simplify and we consider  
  
   \begin{thm} \label{lift.e}{\bf Lifting of endofunctors.} \em 
 For a monad $F$, a comonad $G$ and an endofunctor $T$ on the category $\A$,  
consider the diagrams 
 $$\xymatrix{ 
        \A_F \ar@{-->}[r]^\oT \ar[d]_{U_F} & \A_F \ar[d]^{U_F} \\
      \A \ar[r]^T &\A , }\qquad
    \xymatrix{ 
        \A^G \ar@{-->}[r]^\wT \ar[d]_{U^G} & \A^G \ar[d]^{U^G} \\
      \A \ar[r]^T &\A ,}   
  $$  
where the $U$ denote the forgetful functors.
\end{thm}

 \begin{thm} \Label{lift.mon.mon}  {\bf Monad distributive laws.} \em 
If $T$   
is also a monad, 
a natural transformation $\lambda: FT\to TF$ is called a
{\em monad distributive law}
provided $T$ can be lifted to a monad $\oT$ (in \ref{lift.e}).
These conditions can be formulated by some commutative diagrams 
(e.g.\ \cite[4.4]{WisAlg}). 
In this case $\lambda:FT\to TF$ induces a canonical monad structure on $TF$.
 
  $TF$-modules are equivalent to  $\oT$-modules,
i.e.\ $T$-modules whose structure map is a morphism of $F$-modules. This
means $F$-modules $\alpha:F \arg A\to A$  
as well as $T$-modules $\beta: T\arg A\to A$ inducing commutativity of the
 diagram 
$$\xymatrix{FT\arg A \ar[rr]^{\lambda\narg A} \ar[d]_{F\arg \beta} & &  
TF\arg A\ar[d]^{T\arg \alpha} \\
F\arg A \ar[r]^\alpha & A & \ar[l]_{\beta} T\arg A .} $$
\end{thm}

\begin{thm} \Label{lift.comon.comon}  {\bf Comonad distributive laws.} \em
If  $T$ is a comonad, 
a natural transformation $\varphi: TG\to GT$ is called a
{\em comonad distributive law}
provided $T$ can be lifted to a comonad $\wT$ (in \ref{lift.e}). 
The properties can be expressed by some commutative diagrams 
(e.g.\ \cite[4.9]{WisAlg}). In this case 
    $\varphi:TG\to GT$ induces a canonical comonad structure on $TG$.

  $TG$-comodules are equivalent to $\wT$-comodules, i.e.\ 
objects $A$ that are $G$-comodules $\gamma:A\to G\arg A$ as well as
$T$-comodules $\delta: A\to T\arg A$ inducing commutativity of the 
diagram  
$$\xymatrix{T\arg A\ar[d]_{T\arg \gamma} & \ar[l]_{\delta} A \ar[r]^\gamma &
    G\arg A \ar[d]^{G\arg \delta} \\
TG\arg A \ar[rr]^{\varphi\narg A}  & &  GT\arg A. } $$
\end{thm}

 \begin{thm} \Label{f.entw}  {\bf Mixed distributive laws.} \em 
If $T$ is a comonad, a natural transformation  
$\lambda: FT\to TF$ is called a {\em mixed distributive law} or 
an {\em entwining} provided 
the functor $T$ can be lifted to a comonad on the module category $\A_F$
(equivalently, $F$ can be lifted to a monad on the comodule category $\A^T$). 
Again this can be characterised by some commutative diagrams 
(e.g.\ \cite[5.3]{WisAlg}).

  {\em Mixed bimodules} or {\em $\lambda$-bimodules} are defined as those 
$A\in \Obj(\A)$ with morphisms
$$
\xymatrix{ 
F\arg A\ar[r]^{\quad h} & A \ar[r]^{k\quad} &T\arg A
}
$$ 
such that $(A,h)$ is an $F$-module and $(A,k)$ is a
$T$-comodule satisfying the pentagonal law
$$
\xymatrix{ 
F\arg A \ar[r]^{\quad h} \ar[d]_{F\arg k} & A \ar[r]^{k\quad} &T\arg A \\
FT\arg A \ar[rr]^{\lambda\narg A} && TF\arg A \ar[u]_{T\arg h}. 
} 
$$ 
 
Morphisms between two $\lambda$-bimodules, called {\em bimodule morphisms},
are both $F$-module and $T$-comodule morphisms.
  
These notions yield the category of  $\lambda$-bimodules  which we 
denote by $\A^T_F$.   
It can also be considered as the category of $\overline T$-comodules for the  
 comonad $\overline T: \A_F \to \A_F$ or as the category of 
$\widehat F$-modules  for the monad $\widehat F: \A^T\to \A^T$.  
\end{thm}

In \cite[2.2]{MW} also entwinings of the type $GF \to FG$ are considered, for
a monad $F$ and a comonad $G$. While the compatibility conditions can be
written formally symmetrically to mixed distributive laws in \ref{f.entw},
such entwinings have no interpretation in terms of liftings.

\begin{thm}\Label{dist-adj}{\bf Distributive laws for adjoint functors.} 
 Let $(L,R)$ be an adjoint pair of endofunctors on a category $\A$ with unit
 $\eta$ and counit $\adjcu$, 
and $F$ be an endofunctor on $\A$. Consider a natural transformation
$\psi:LF\to FL$ and set 
$$\tilde \psi=  RF\adjcu  \circ R\psi R   \circ \eta FR  : FR\to RF.$$
\begin{zlist}
\item  If $L$  and $F$ are monads, then $\psi$ is a monad
  distributive law if and only if $\tilde \psi$ is an entwining (mixed
  distributive law). 
\item  If $L$ is a monad and $F$ a comonad, then $\psi$ is an entwining
  (mixed distributive law) if and only if $\tilde \psi$ is a comonad
  distributive law. 
\item  If $L$ is a comonad and $F$   a monad, then $\psi$ is an entwining if
  and only if $\tilde \psi$ is a monad distributive law.
\item  If $L$ and $F$  are comonads, then $\psi$ is a comonad
  distributive law if and only if $\tilde \psi$ is an entwining. 
\end{zlist}  
\end{thm}
\begin{proof}  
All these claims are easily checked by using that the structure maps of the
adjoint monad-comonad (or comonad-monad) pair $(L,R)$ are mates under
adjunctions, together with naturality and the triangle identities. 
Details are left to the reader. 
\end{proof}
 
Combining the correspondences in \ref{dist-adj}(1) and (2) with the
isomorphism of module and comodule categories in \ref{adj-mon}(1), further
isomorphisms, between categories of mixed bimodules, can be derived. 
Note, however, that since the entwinings occurring in parts (3) and (4) of
\ref{dist-adj} can not be translated to liftings, these claims do not lead to
similar conclusions.
 
\begin{thm}\Label{dist-mod}{\bf Modules and distributive laws.} 
Let $L$ be a monad with right adjoint comonad $R$ on a category $\A$.
\begin{zlist}
\item Let $G$ be a comonad with a mixed distributive law $\lambda: LG\to GL$.
Then the category of mixed $(L,G)$-bimodules $\A_L^{G}$ 
is isomorphic to the 
category of $(R,G)$-bicomodules $\A^{GR}$ 
 (see \ref{lift.comon.comon})
defined by the associated comonad distributive law $\tilde\lambda: GR \to RG$
(see \ref{dist-adj}(2)).

\item Let $F$ be a monad with a mixed distributive law $\tau: FR\to RF$.
Then the category of mixed $(F,R)$-bimodules $\A_F^R$
is isomorphic to the category of $(L,F)$-bimodules  
 $\A_{FL}$ 
defined by the monad distributive law  $\tilde\tau: LF\to FL$
(see \ref{dist-adj}(1)).
\end{zlist}
\end{thm}

\begin{proof} 
(1)  The mixed distributive law $\lambda:LG \to  GL$ yields a lifting of $G$ to
  a comonad  ${\overline G}$ on the category $\A_L$ of $L$-modules. Moreover,
  it determines a comonad distributive law $\tilde\lambda:GR \to RG$  
 which is equivalent to a lifting of $G$ to a comonad ${\widehat G}$ 
 on the category $\A^R$ of $R$-comodules. By \ref{adj-mon}(1), $\A_L$ and
 $\A^R $ are isomorphic, and this isomorphism obviously `intertwines' the 
comonads ${\widehat G}$ and ${\overline G}$. Thus the isomorphism 
$\A_L \iso \A^R$ lifts to an isomorphism between the categories of ${\widehat
  G}$-comodules and ${\overline G}$-comodules. By characterisation of
${\widehat G}$-comodules as comodules for the composite comonad $GR$, and
characterisation of ${\overline G}$-comodules as mixed $(L,G)$-bimodules, we
obtain the isomorphism claimed.

In fact, it can also be computed directly that a mixed bimodule 
$$
\xymatrix{ 
L\arg A\ar[r]^{\quad h} \ar[d]_{L\arg k} & A \ar[r]^{k\quad} & G\arg A \\
LG\arg A \ar[rr]^{\lambda\narg A} && GL\arg A \ar[u]_{G\arg h}, 
} 
$$ 
transforms to an $(R,G)$-bicomodule  with commutative diagram
$$
\xymatrix{ 
R\arg A\ar[d]_{R\arg k} & \ar[l]_{\quad h'}   
     A \ar[r]^{k\quad} & G\arg A \ar[d]^{G\arg h'} \\
   RG\arg A & & \ar[ll]^{\tilde\lambda\narg A}   GR \arg A  , } 
$$ 
where $h'$ is the mate of $h$ under the adjunction $(L,R)$, cf. \ref{adj-mon}.

(2) is shown similarly to (1).
\end{proof}

\begin{thm}\Label{F-act}{\bf $\nF$-actions on functors.} \em
Let $\A$ and $\B$ be categories.
Given a monad $\bF=(F,m,i)$ on $\A$, any functor $R : \B \to \A$
is called a
{\em (left) $\nF$-module} if there exists a natural transformation
$\varphi: FR \to R$ satisfying associativity and 
unitality conditions (corresponding to those required 
for objects, see \cite[3.1]{MW}).
Clearly, for any functor $R:{\mathbb B}\to {\mathbb A}$,  
$(FR,  m R)$ is an $\nF$-module functor.
\end{thm}

\begin{thm}\Label{F-Gal}{\bf $\nF$-Galois functors.} \em
For a monad $\bF$ on a category $\A$ and any functor $R : \B \to \A$,
consider the diagram 
 $$\xymatrix{ 
  \A_\nF \ar@{-->}[r]^{\overline{L}} 
 & \B \ar@{-->}[r]^{\overline{R}}  
\ar@{=}[d]
& \A_\nF \ar[d]^{U_\nF} \\
\A \ar@{-->}[r]^L \ar@{-->}[u]^{\phi_\nF}& \B \ar[r]^R  & \A . }
  $$    
As a particular instance of \ref{Apple}, 
there exists some functor $\overline{R}$ making the right square commutative
if and only if $R$ has an $\nF$-module structure $\varphi:FR\to R$ (see
\ref{F-act}). 

If $R$ has a left adjoint $L:\A\to \B$, with unit $\eta$ and counit $\adjcu$ of
the adjunction, then there is a monad morphism 
$$\can: F\stackrel{F\eta}\lra FRL\stackrel{\varphi L}\lra RL.$$ 

We call an $\nF$-module functor $R$ an {\em $\nF$-Galois functor} if it has a
left adjoint and $\can$ is an isomorphism.  
\medskip

Consider an $\nF$-module functor $R:\B \to \A$ with $F$-action $\varphi$, a left
adjoint $L$, unit $\eta$ and counit $\adjcu$ of the adjunction. If
$\B$ admits coequalisers of the parallel morphisms $L\arg \rro_X$ and
$\adjcu\narg{ L\arg X}\circ L\arg{\varphi \narg {L\arg X}} \circ LF\arg{\eta
  \narg X}:LF\arg X\to L\arg X$, for any object $(X,\rro_X)$ in $\A_\nF$,
then 
this coequaliser yields the left adjoint $\overline L(X,\rro_X)$ of
${\overline R}$ (see left square). By uniqueness of the adjoint,
$\overline{L}\phi_\nF \iso L$ (see \ref{Apple}).   
Denoting the coequaliser natural epimorphism $LU_\nF \to {\overline L}$ by $p$,
the unit of the adjunction $({\overline L},{\overline R})$ is the unique natural
morphism 
${\overline \eta}:\id_{{\mathbb A}_\bF} \to {\overline R}\,{\overline L}$ such
that $U_\nF {\overline \eta}= Rp \circ \eta U_\nF$. The counit is the unique
natural morphism ${\overline \adjcu}:{\overline L}\, {\overline R}\to
\id_{\mathbb B}$, such that ${\overline \adjcu}\circ p {\overline R}  =
\adjcu$. 
\smallskip
 
If $\B$ has coequalisers of all parallel morphisms,  
then the following are equivalent (dual to  \cite[Theorem  3.15]{MW}):
\begin{blist}
\item $R$ is an $\nF$-Galois functor;
\item the unit of $(\overline{L},\overline{R})$ is an isomorphism for 
   \begin{rlist}
     \item all free $\nF$-modules (i.e.\ modules in the Kleisli category of
       $\bF$), or  
     \item all $U_\nF$-projective $\nF$-modules.
   \end{rlist} 
\end{blist} 
\end{thm}

 {}From  \cite{Dub:Kan} and \cite[3.3 Theorem 10]{TTT} 
 we recall a result of central importance in our setting
in the form it can be found in \cite[Theorem 1.7]{Pepe:Comonad&coring}.

\begin{thm}\Label{Beck}{\bf Beck's theorem.}  
Consider a monad $F$ on a category $\A$ and an $F$-module functor $R:\B \to
\A$. Then the induced lifting ${\overline R}: \B \to \A_F$ in
\ref{F-Gal} is an equivalence if and only if the following hold:  
\begin{rlist}
\item $R$ is an $F$-Galois functor,
\item $R$ reflects isomorphisms,
\item ${\mathbb B}$ has coequalisers of $R$-contractible coequaliser pairs
  and $R$ preserves them.
\end{rlist}
\end{thm}
\medskip

The Galois property also transfers to adjoint functors.

\begin{thm}\Label{quest}{\bf Proposition.}   
Consider an adjoint pair $(F,G)$ of endofunctors on a category $\A$.
Let $T:\B \to \A$ be a functor which has both a left adjoint $L$ and a
right adjoint $R$. 
\begin{zlist}
\item If $F$ is a monad (equivalently, $G$ is a comonad) then 
$T$ is an $F$-Galois functor as in \ref{F-Gal} if and
only if it is a $G$-Galois functor (in the sense of \cite[Definition
3.5]{MW}). 

\item If $F$ is a comonad (equivalently, $G$ is a monad) then 
$R$ is a $G$-Galois functor as in \ref{F-Gal} if and
only if $L$ is an $F$-Galois functor (in the sense of \cite[Definition
3.5]{MW}). 
\end{zlist}
\end{thm}

\begin{proof}
Denote the unit of the adjunction $(F,G)$ by $\eta$ and its counit by
$\adjcu$. Denote furthermore the unit and counit of the adjunction $(L,T)$ by
$\eta_L$ and $\adjcu_L$, respectively, and for the unit and counit of the
adjunction $(T,R)$ write $\eta_R$ and $\adjcu_R$, respectively.

(1) A bijective correspondence between $F$-actions $\varphi_T$ and
  $G$-coactions $\varphi^T$ on $T$ is given by $\varphi^T:=G \varphi_T \circ
  \eta T$.
The canonical comonad morphism corresponding to $\varphi^T$ comes out as
$$\xymatrix{
{\widetilde{\can}}\ :\ 
TR \ar[r]^-{\eta T R} & 
G F T R \ar[r]^-{G \varphi_T R} &
G TR \ar[r]^-{G {\adjcu}_R} & G\ .
}$$
Comparing it with the canonical monad morphism $\can:F \to TL$ in 
\ref{F-Gal}, they are easily seen to be mates under the adjunctions $(F,G)$
and
$(TL,TR)$. That is, 
$$
{\widetilde{\can}}= G {\adjcu}_R  \circ  G T \adjcu_L
{R}   \circ  G \can T{R}  \circ  \eta T {R}. $$
Thus ${\widetilde{\can}}$ is an isomorphism if and only if
$\can$ is an isomorphism.

(2) A bijective correspondence between $G$-actions $\varphi_R:GR\to R$ and
$F$-coactions $\varphi^L:L\to FL$ is given by
$$
\varphi^L:= {\adjcu_L FL}\circ {L\adjcu_R TFL} \circ {LG \varphi_R TFL} \circ
{LTG \eta_R FL} \circ {LT \eta L} \circ {L\eta_L}.
$$
The canonical comonad morphism ${\widetilde\can}:LT \to F$ corresponding to
$\varphi^L$, and the canonical monad morphism $\can:G \to RT$ corresponding to
$\varphi_R$, turn out to be mates under the adjunctions $(LT,RT)$ and
$(F,G)$. That is, 
$$
{\widetilde \can}=\adjcu_L F \circ L\adjcu_R TF \circ LT \can F \circ LT \eta.
$$
This proves that $\can$ is a natural isomorphism if and only if ${\widetilde
  \can}$ is a natural isomorphism, as stated.
\end{proof}

\section{Rings in module categories}

 Let $A$ be an associative ring with unit.
In this section we study the relationship between ring extensions of $A$ and
monads on the category $\M_A$ of right $A$-modules. 

 \begin{thm}\Label{ring.ex}{\bf $A$-rings.} \em
A ring $B$ is said to be an {\em $A$-ring} provided there
is a ring morphism $\unit:A\to B$. This is equivalent to saying that
 $B$ is an  $\AA $-bimodule with $A$-bilinear multiplication and unit,
$$
\product:B\ot_A B \to B, \quad \unit:A\to B,$$
 inducing commutative diagrams for associativity and unitality.
 
A right {\em $B$-module} is a right $A$-module $M$ with an $A$-linear map 
  $$\rro_M: M\ot_A B \to M,\quad m\ot b\mapsto mb,$$
satisfying the associativity and unitality conditions.
  $B$-module morphisms $f:M\to N$  are $A$-linear maps and
 $f\circ\rro_M  = \rro_N\circ (f \ot_A\id_B)$.
 
The category of right $B$-modules is denoted by $\M_B$. It is isomorphic
to the module category over the ring $B$ and thus is  
an abelian category with $B$ as a projective generator. 
\end{thm}

\begin{thm}\Label{adj-BHom}{\bf Adjointness of $-\ot_A B$ and $\Hom_A(B,-)$.}
  \em 
As an endofunctor on $\M_A$, $-\ot_AB$ is left adjoint to the 
endofunctor $\Hom_A(B,-)$ with unit and  counit  
$$\begin{array}{rl} 
\eta_X: X \to \Hom_A(B, X\ot_A B), & x\mapsto [b\mapsto x\ot b],\\[+1mm] 
\adjcu_N: \Hom_A(B,N)\ot_A B \to N, & f\ot b \mapsto f(b).
\end{array}$$ 
\end{thm}

\begin{thm}{\bf Remark.} \Label{rem:gen}\em
Since $A$ is a generator in $\M_A$ and
$-\ot_AB$ preserves direct sums and epimorphisms, 
the functor $-\ot_AB: {\M_A} \to {\M_A}$ is fully determined by the value at
$A$, that is by $A\ot_A B\iso B$. 
Similarly, a natural
transformation $\varphi: -\ot_A B \to -\ot_A B'$ between such `tensor
functors' is equal to $-\ot_A \varphi_A$, where $\varphi_A:B \to B'$ is an
$A$-bimodule map.

In general,  $\Hom_A(B,-)$ is not
determined by ${B^*}=\Hom_A(B,A)$, unless it preserves direct sums and
epimorphisms, that is, unless $B$ is a finitely generated and projective right
$A$-module.  
However, 
 $\Hom_A(B,-)$ is determined by  $\Hom_A(B,Q)$ for any cogenerator
 $Q \in \M_A$ since it is left exact and preserves direct products.
For a natural transformation $\varphi:\Hom_A(B,-)\to \Hom_A(B',-)$
between `Hom functors', it follows by the Yoneda Lemma that
$\varphi=\Hom_A(\varphi_B(\id_B),-)$, where $\varphi_B(\id_B):B'\to B$ is an
$\AA $-bimodule map. 
\end{thm}

\begin{thm}\Label{mon-com}{\bf Monad-comonad.} 
 For an $\AA $-bimodule $B$ the following are equivalent:
\begin{blist}
 \item $(B,\product ,\unit )$ is an $A$-ring;
 \item $-\ot_AB : \M_A \to \M_A$ is a monad;
\item  $\Hom_A(B,-): \M_A \to \M_A$ is a comonad.
\end{blist}
\end{thm}

 \begin{proof} 
(a)$\LRa$(b).  An $A$-ring $(B,\product ,\unit )$ determines a monad $(-\ot_A
   B,-\ot_A \product , 
   -\ot_A \unit )$ on $\M_A$ and this is a bijective correspondence in light of
   Remark \ref{rem:gen}. The equivalence (b)$\LRa$(c) follows by \ref{adj-mon}. 
 
Explicitly, for an $A$-ring $(B,\product ,\unit )$,
$\Hom_A(B,-)$ is a comonad by coproduct and counit
$$
\begin{array}{c} 
\Hom_A(B,-)\stackrel{\Hom_A(\product ,-)}\lra   
\Hom_A(B\ot_A B,-)\stackrel{\iso}\lra   \Hom_A(B,\Hom_A(B,-)), \\[+1mm] 
\Hom_A(B,-) \stackrel{\Hom_A(\unit ,-)}\lra \Hom_A(A,-)\stackrel\iso
\lra\id_{\M_A}.  
\end{array} 
$$
\end{proof}

Adjointness of the free and forgetful functors for the monad $-\ota B$ is 
just the isomorphism 
$$ \Hom_B(-\ot_A B,N) \to \Hom_A(-,N), \quad f\mapsto   f \circ (- \ot_A
\unit). $$  

We write $[B,-]=\Hom_A(B,-)$ for short. 
Comultiplication and counit of the comonad $[B,-]$ in \ref{mon-com}(c) are
denoted by $[\product,-]$ and $[\unit,-]$, respectively.

The following was pointed out in \cite[page 141]{TTT}. 

 \begin{thm}\Label{[B].mod}{\bf $B$-modules are $\Hom_A(B,-)$-comodules.}
For any $A$-ring $B$,
the category of right $B$-modules is isomorphic to the category of
$\Hom_A(B,-)$-comodules, that is, there exists an isomorphism
 $$ {\M_B}  \stackrel\iso\lra {\M^{[B,-]}}.$$
\end{thm}
\begin{proof} 
This is a special case of \ref{adj-mon}(1). 
Here the isomorphism has the following form.
 Given a $B$-module
 $\rro_N:N\ot_A B \to N$, applying $\Hom_A(B,-)$ and composing with the 
unit of the adjunction yields 
  $$\xymatrix{
N\ar[r]^-{\eta_N} & \Hom_A(B,N\ot_AB)
 \ar[rr]^-{\Hom(B,\rro_N)}& & \Hom_A(B,N),}$$ 
imposing a $\Hom_A(B,-)$-comodule structure on $N$.
 
Conversely, a comodule structure map
  $\rro^N:N \to\Hom_A(B,N) $ induces
 $$\xymatrix{ 
N\ot_A B \ar[rr]^-{\rro^N \ot_A B} && \Hom_A(B, N) \ot_A B 
\ar[r]^-{\vareps_N} & N,}$$
defining a $B$-module structure on $N$.
\end{proof}

\section{Corings in module categories}

Again $A$ denotes an associative ring with unit. 
To any $A$-coring $\cC$ we can associate a comonad $-\ot_A \cC$ and a monad
$\Hom_A(\cC,-)$ on the category $\M_A$ of right $A$-modules. Here we 
consider the relationship between $-\ot_A \cC$-comodules (i.e.\ {\em
  $\cC$-comodules}) and $\Hom_A(\cC,-)$-modules (i.e.\ {\em
  $\cC$-contramodules}).  
 
\begin{thm}{\bf $A$-Corings.} \em
An $A$-coring is an $\AA $-bimodule $\C$ with $A$-bilinear maps
 $$\cop: \C\to \C\ot_A \C, \quad \counit:\C\to A$$  
satisfying coassociativity and counitality conditions. 
\end{thm} 

Similar to the characterisation of $A$-rings in \ref{mon-com}, 
 \ref{adj-mon} implies the following characterisation of corings:

\begin{thm}\Label{coring}{\bf Corings.}
For an $\AA $-bimodule $\C$, the following are equivalent:
\begin{blist}
\item $(\C,\cop,\counit)$ is an $A$-coring;
\item $- \ot_A \C:{\M}_A\to {\M}_A$ induces a comonad;
\item $\Hom_A(\C,-):{\M}_A\to {\M}_A$ induces a monad.  
\end{blist} 
\em Writing $\Hom_A(\C,-)=[\C,-]$,
 the related monad is $([\C,-], [\cop,-], [\counit,-])$. 
\end{thm}

In the rest of this section $\C$ will be an $A$-coring. 
We first recall properties of the category of comodules for the 
related comonad (see \cite{BW}). 

 \begin{thm}\Label{cat.C}{\bf The category $\M^\C$.}
The comodules for the comonad $- \ot_A \C:\M_A\to \M_A$
are called {\em right $\C$-comodules} and their category is denoted by
${\M^\C}$.
\begin{zlist}
\item  $\M^\C$ is an additive category with coproducts and cokernels.
\item The (co)free functor $- \ot_A \C$ is right adjoint to the forgetful 
functor by the isomorphism 
  $$ \Hom^\C(M,X\ot_A \C) \to \Hom_A(M,X), \quad f\mapsto (\id_X
\ot_A\counit)\circ f.$$  
\item For any generator  $P\in \M_A$, $P \ot_A \C$ is a subgenerator in  $M^\C$,
  in particular, $\C$ is a subgenerator in  $\M^\C$.
\item For any cogenerator  $Q\in \M_A$, $Q\ot_A \C$ is a cogenerator in
  $\M^\C$.      
\item For any injective  $X\in \M_A$, $X \ot_A \C$ is injective in $\M^\C$. 
\item For any monomorphism $f:X\to Y$ in $\M_A$, $f \ot_A \id_\C: X\ot_A \C
  \to Y \ot_A \C$ is a monomorphism in $\M^\C$. 
\item For any family $X_\lambda$ of $A$-modules, 
     $\prod_\Lambda X_\lambda\ot_A \C$ 
     is the product of the $\C$-comodules  $X_\lambda\ot_A \C$.
\item  $\C$ is a flat left $A$-module if and only if  
    monomorphisms in $\M^\C$ are injective maps.
\end{zlist}
\end{thm}

Left comodules of an $A$-coring $\C$ are defined symmetrically to the right
comodules in \ref{cat.C}, as comodules of the comonad $\C\ot_A -$ on the
category of left $A$-modules. 
Furthermore, if $\C$ is an $A$-coring and $\D$ is a $B$-coring, then 
we can consider the (composite) comonad 
$\C \ot_A - \ot_B \D$ on the category of $(A,B)$-bimodules. 
Its comodules are called {\em $(\C,\D)$-bicomodules}. Equivalently, a
$(\C,\D)$-bicomodule is a left $\C$-comodule and a right $\D$-comodule  
such that the right $\D$-coaction is a left $\C$-comodule map,
cf. \ref{lift.comon.comon}. 
The category of $(\C,\D)$-bicomodules is denoted by ${}^\C\M^\D$.

While $\C$-comodules are well studied in the literature (e.g.\ \cite{BW}),
 $\kC$-modules have not attracted so much attention so far.  
They were addressed by Eilenberg-Moore in \cite{EM} as
 {\em $\C$-contramodules} and reconsidered 
recently by Positselski \cite{Pos} in the context of 
semi-infinite cohomology.

\begin{thm}\Label{MC}{\bf The category $\M_\kC$.}   
The modules for the  monad $\kC:\M_A \to \M_A$ 
are right $A$-modules $N$ with some
$A$-linear map $[\C,N]\to N$ subject to associativity and unitality
conditions. Their category is denoted by ${\M_\kC}$.
\begin{zlist}
\item  $\M_\kC$ is an additive category with products and kernels.
\item The (free) functor $\kC: {\M_A}\to {\M_\kC}$ is left adjoint to the 
      forgetful functor by the isomorphism   
$$\Hom_\kC([\C,X],M) \to \Hom_A(X,M), \quad f\mapsto f \circ [\counit,X].$$ 
\item For any generator  $P\in M_A$, $[\C,P]$ is a generator in  $\M_\kC$;
  in particular, $\C^*=[\C,A]$ is a generator in $\M_\kC$ and 
  $\Hom_\kC ([\C, A],M) \iso M$.
\item For any cogenerator  $Q\in \M_A$, $[\C,Q]$ is 
    a weak subgenerator in $\M_\kC$, that is, every object of $\M_\kC$
    is a subfactor of some product of copies of $[\C,Q]$ (compare
    \cite[1.6]{WisCot}). 
\item For any projective  $Y\in \M_A$, $[\C,Y]$ is projective in $\M_\kC$. 
 \item For any epimorphism $h:X\to Y$ in $\M_A$,
       $[\C,h]: [\C,X]\to [\C,Y]$
      is an epimorphism (not necessarily surjective) in $ \M_\kC$.
\item For any family $X_\lambda$ of $A$-modules, 
     $[\C,\bigoplus_\Lambda X_\lambda]$ 
   is the coproduct of the $[\C,-]$-modules 
    $[\C, X_\lambda]$.
\item $\C$ is a projective right $A$-module if and only if 
     epimorphisms in $\M_\kC$ are surjective maps.
\end{zlist}
\end{thm}

\begin{proof} The proof is similar to the comodule case. 
Some of the assertions can also be found in \cite{Pos}.
\end{proof} 
 
\begin{thm}{\bf Right and left contramodules.} \Label{rem:left}\em
In \ref{MC}, modules of the monad $[\C,-]\equiv \Hom_{-,A}(\C,-)$ on the
category $\M_A$ of right $A$-modules are considered. Symmetrically, an
$A$-coring $\C$ determines a monad $\Hom_{A,-}(\C,-)$ also on the category
${}_A \M$ of left $A$-modules. In situations when both monads on the
categories of right and left $A$-modules occur at the same time, we use the
following terminology to distinguish between their modules.
Modules for the monad $[\C,-]\equiv \Hom_{-,A}(\C,-)$ on $\M_A$ are called {\em
  right} $\C$-contramodules, and modules for the monad $\Hom_{A,-}(\C,-)$ on
${}_A\M$ are called {\em left} $\C$-contramodules. 
If not specified otherwise, we mean by contramodules {\em right}
contramodules, throughout.
\end{thm}

We saw in \ref{[B].mod} that for any $A$-ring $B$, the categories $\M_B$ and
$\M^{[B,-]}$ are isomorphic (see also \ref{adj-mon}(1)). 
In view of 
the asymmetry of assertions (1) and (2) in \ref{adj-mon}, 
the corresponding 
statement for corings is no longer true and we will come back to this 
question in \ref{corresp}. 
So far we know from \ref{adj-mon}(2):

\begin{thm}\Label{Klei}{\bf Related Kleisli categories.} 
For any $A$-coring $\C$, the Kleisli categories of $- \ot_A \C$
and $\kC$ are isomorphic by the isomorphisms for $X,Y\in \M_A$,
$$\begin{array}{rcl}
\Hom^\C(X \ot_A \C,Y\ot_A \C)&\iso& \Hom_A(X \ot_A \C, Y) \\
   &\iso & \Hom_A(A,[\C,Y]) \\
 &\iso & \Hom_\kC([\C,X],[\C,Y]).
\end{array}$$
\end{thm} 

 Recall that for any $A$-coring $\C$, the right dual $\C^*=\Hom_A(\C,A)$ 
has a ring structure by the convolution product for $f,g\in {\C}^*$,
$f*g= f \circ(g \ot_A \id_\C) \circ\cop$
(convention opposite to \cite[17.8]{BW}).
Similarly a product is defined for the left dual ${^*\C}$.

The relation between  $\C$-comodules and modules over the dual ring of $\cC$ is 
well studied (see e.g.\ \cite[Section 19]{BW}). 

\begin{thm}\Label{com.*}{\bf The comonads $- \ot_A \C$ and $[{}^*\C,-]$.} 
The comonad morphism  
$$\alpha :- \ot_A \C\to \Hom_A({}^*\C,-), \quad - \ot c\mapsto 
[f\mapsto -f(c)],$$ 
yields a faithful functor
$$ \begin{array}{rcl}
G_\alpha:\qquad {\M^\C}  &\lra & \M^{[{}^*\C,-]}\iso \M_{{}^*\C} \ ,\\
\big(N\stackrel{\rro^N}{\lra} N \ot_A \C\big) &\longmapsto &
\big(N\stackrel{\rro^N}{\lra} N\ot_A \C \stackrel{\alpha_N}\lra
\Hom_A({}^*\C,N)\big) 
\end{array}$$ 
 and the following are equivalent: 
\begin{blist}
\item $\alpha_N$ is injective for each $N\in \M_A$;
\item $G_\alpha$ is a full functor;
\item $\C$ is a locally projective left $A$-module.
\end{blist}  

If these conditions are satisfied, $\M^\C$ is equal to $\sigma[\C_{{}^*\C}]$,
the full subcategory of $\M_{{}^*\C}$ subgenerated by $\C$. 
\end{thm}
 
Similar to \ref{com.*},  $\C$-contramodules
can be related to $\C^*$-modules.

\begin{thm}\Label{monC*}{\bf The monads $\kC$ and $-\ot_A\C^*$}.  
The monad morphism
 $$
\beta : - \ot_A \C^* \to \Hom_A(\C,-), \quad -\ot f\mapsto [c\mapsto -f(c)], $$ 
 yields a faithful functor
$$ \begin{array}{rcl}
F_\beta:\qquad {\M_\kC} &\lra & \M_{\C^*}\ ,\\
\big(\Hom_A(\C,M) \stackrel{\rro_M}{\lra} M\big) &\mapsto & 
\big(M\ot_A \C^* \stackrel{\beta_M}\lra \Hom_A(\C,M)\stackrel{\rro_M}{\lra} M
\big),
  \end{array}$$
associating to a $\kC$-module $M$ the same object with 
a ${\C^*}$-module structure.
 The following are equivalent: 
\begin{blist}
\item $\beta$ is surjective for all $M\in \M_A$;
\item $F_\beta$ is a full functor;
\item $\C$ is a finitely generated and projective right $A$-module;
\item $F_\beta$ is an isomorphism.
\end{blist}  
\end{thm}
 
In general $\C$ is not a $\kC$-module and $[\C,A]$ is not a $\C$-comodule. 
In fact,  $[\C,A]\in {}^\C\M$ holds provided $\C$ is finitely generated and
projective as a right $A$-module.

\section{Functors between co- and contramodules}

Categories of comodules and contramodules have complementary
features. Therefore, it is of interest to find $A$-corings $\C$ and 
$B$-corings $\D$ (over possibly different base rings) such
that the category of $\D$-comodules and that of $\kC$-modules are
equivalent. 
As we will see  in \ref{thm:EM}, functors between these categories are 
provided by bicomodules. It  turns out that the question, when they provide an
equivalence, fits the standard problem in (categorical) descent theory.

Since comodules for the trivial $B$-coring $B$ are simply $B$-modules, 
our considerations include the particular case when the 
category of $\kC$-modules is equivalent to the category of $B$-modules. 
Dually, when the coring $\C$ is trivial (i.e.\ equal to $A$), the problem
reduces to a study of equivalences between $A$-module and $\D$-comodule
categories. This question is already discussed in the literature, see
e.g.\ \cite{Joost:Equi}, \cite{Pepe:Comonad&coring}. 

Throughout this section $\C$ is an $A$-coring and $\D$ a $B$-coring
for rings $A$ and $B$. The following observation was made in 
\cite[5.1.2]{Pos}.

\begin{thm}\Label{thm.com.contra}
{\bf $\kC$-modules induced by $\C$-comodules.}
Let $N$ be a $(\cC,\cD)$-bico\-module with left $\cC$-coaction $\Nro$.
For any  $Q\in \M^\cD$, there is an isomorphism  
$$\varphi:\Hom_A(\C, \Hom^\cD(N,Q)) \to \Hom^\cD(\C\ot_A N, Q), \;
   h\mapsto [ c\ot m \mapsto h(c)(m)],$$ 
(see e.g.\ \cite[18.11]{BW}). 
Then the right $A$-module $N^Q := \Rhom \cD N Q$ is a $\kC$-module by
$\alpha _{N^Q}$:
$$
\xymatrix{
 \Hom_A(\C, \Hom^\cD(N,Q)) \ar[r]^-{\varphi}&
\Hom^\cD(\C\ot_A N, Q) \ar[rr]^-{\Hom^\cD(\Nro,Q)}&&\Hom^\cD(N, Q) . 
}
$$
Thus there is a bifunctor 
$\Rhom \cD - - : ({}^\cC\M^\cD)^{op}\times \M^\cD \to \M_\kC$,  
$$(N,Q)\mapsto (N^Q, \alpha_{N^Q}), \quad (f,g)\mapsto \Rhom \cD f g.$$
\end{thm}

\begin{proof} The identification of $\rhoma  \cC {\Rhom \cD N Q}$ with 
$\Rhom \cD {\cC\ota N} Q$, yields $\alpha_{N^Q} = \Rhom \cD \Nro Q$.
Since the left $\cC$-coaction $\Nro$ of a $(\cC,\cD)$-bicomodule $N$ is right
$\cD$-colinear, $\alpha_{N^Q}(f)$ is right $\cD$-colinear,
for all $f\in \Rhom \cD {\cC\ota N} Q$.
 Hence $\alpha_{N^Q}$ is well-defined. The left
$A$-linearity of $\Nro$ implies that $\alpha_{N^Q}$ is right $A$-linear.  
  By coassociativity of $\Nro$, one 
concludes on associativity of $\alpha_{N^Q}$.
Similarly, using the counitality of $\Nro$, one proves the commutativity of
the triangle diagram in \ref{monad} for $\alpha_{N^Q}$. Again a similar
computation  
yields that for a morphism $f:N\to N'$ of $(\cC,\cD)$-bicomodules, 
$\Rhom \cD f Q$ is a morphism of $\kC$-modules. 
If $g: Q\to P$ is a morphism of right $\cD$-comodules, 
then $\Rhom \cD N g$ is a morphism of $\kC$-modules since  
 $\alpha_{N^Q} = \Rhom \cD \Nro Q$ implies that $\alpha_{N^Q}$ 
 is  natural in $Q$. 
\end{proof}

In a symmetric way, for any left $\cD$-comodule $Q$ and a
$(\cD,\cC)$-bicomodule $N$ (with $\cC$-coaction $\roN: N\to N\ota \cC$),
$\Hom^{\cD} (N, Q)$ is a left $\C$-contramodule 
by $\Hom^{\D}(\roN,Q)$.   

If $N$ is just a left $\C$-comodule we tacitly assume $\D=B={\End^\C}(N)$
to apply the preceding notions and results. 

\begin{thm}\Label{HomC}{\bf Corollary.}
\begin{zlist}
\item Let $N$ be a left $\C$-comodule with $B=\End^\C(N)$.
For any subring $B'\subset B$ and $Q\in \M_{B'}$,
$\Hom_{B'}(N,Q)$ is a $\kC$-module. 
\item For any $Q\in \M^\C$, $\Hom^\C(\C,Q)$ is a $\kC$-module.
\end{zlist}
\end{thm}

\begin{thm}\Label{contra.prod}{\bf Contratensor product.} \em
For any $(\C,\D)$-bicomodule $N$, 
the construction in \ref{thm.com.contra} yields a functor 
$${\oHom^\cD}(N,-): \M^\cD \to \M_\kC,$$
 inducing the commutative diagram of (right adjoint) functors 
$$ 
\xymatrix{
 \M^\D \ar[rr]^-{\oHom^\D(N,-)} 
\ar@{=}[d]
&&\M_\kC \ar[d]^-{U_\kC}\\
\M^\D \ar[rr]_-{\Hom^\D(N,-)}&& \M_A\ .}
$$ 
Since $\Hom^\D(N,-)$ has the left adjoint $-\ot_AN$ and
$\M^\D$ has coequalisers, 
it follows from \ref{F-Gal} that ${\oHom}^\cD (N,-)$
also has a left adjoint which comes out as follows   
 (see \cite{Pos}). 

 For any $(\cC,\cD)$-bicomodule  $(N,\Nro,\roN)$ and 
$\kC$-module  $(M,\alpha_M)$, the {\em contratensor product}, $M\otc N$ is
 defined as the coequaliser  
$$
\xymatrix{ \rhoma  \cC M\ota N  \ar@<0.5ex>[rr] \ar@<-0.5ex>[rr] && M\ota
  N\ar[rr] && M\otc N ,} 
$$
where the coequalised maps are $f\ot n\mapsto (f\ota \id_N) \circ\Nro (n)$ and
$\alpha_M\ota \id_N$.  
Projection of an element $m\ot n$ to $M\otc N$ is denoted by $m\otc n$. 

As a coequaliser of right $\cD$-comodule maps, $M\otc N$ is a right
$\cD$-comodule, and thus defines a functor $-\otc N : \M_\kC\to \M^\cD$.   
Note that this coequaliser splits in $\M^\cD$ provided $(M,\alpha_M)$ is
$U_\kC$-projective.
\end{thm}

\begin{thm}\Label{thm:EM}{\bf Functors between comodules and contramodules.}
  Any $(\cC,\cD)$-bicomodule $N$ induces an adjoint pair of functors
$$
-\ot_\kC N: \M_\kC \to \M^\D, \quad \Hom^\D(N,-):\M^\D \to \M_\kC, 
$$
that is, for $M\in \M_\kC$ and $P\in \M^\cD$, there is an isomorphism
 $$
 \Rhom \cD {M\otc N} P \iso \rhom \kC M {\Rhom \cD N P}.
 $$
 
 Conversely, any right adjoint functor $F:\M^\D \to
\M_\kC$ is naturally isomorphic to $\Hom^\D(N,-)$, for an
appropriate $(\C,\D)$-bicomodule $N$. 
\end{thm}

\begin{proof}
In view of the discussion in  \ref{F-Gal}, 
for $\kC$-modules 
$M$, the unit of the adjunction is given by,
$$
\eta_M
: M\to \Rhom \cD N {M\otc N}, \; m\mapsto [n\mapsto m\otc n].
$$
Also by  \ref{F-Gal}, the counit of the adjunction comes out (and is in
particular well defined) as
$$
{\adjcu}_Q: \Rhom\cD N Q\otc N\to Q, \quad f\otc n\mapsto f(n),
$$
for all right $\cD$-comodules $Q$.

Conversely, assume that 
$F:\M^\D \to \M_\kC$ has a left adjoint. Then so does the
composite $F':=U_\kC\circ F:\M^\D\to \M_A$, in light of
\ref{MC}. 
Hence it follows by \cite[Theorem 3.2]{Joost:Equi} that there exists an
$(A,\D)$-bicomodule $N$ such that $F'$ is naturally isomorphic to
$\Hom^\D(N,-)$. Moreover, by construction, for any $Q\in {\M}^\D$, 
$\Hom^\D(N,Q)$ is a $\kC$-module via some action
$\kappa_Q: \Hom_A(\C, \Hom^\D(N,Q)) \to \Hom^\D(N,Q)$,
and for $q\in \Hom^\D(Q,Q')$, $\Hom^\D(N,q)$ is a morphism of
$\kC$-modules. This amounts to saying that $\kappa_{(-)}$ is a natural
transformation $\Hom^\D(\C \otimes_A N,-)\to\Hom^\D(N,-)$. 
Therefore, it follows by the Yoneda Lemma that there is
an $(A,\D)$-bicomodule map $\tau:N \to \C \otimes_A N$, such that
$\kappa_Q=\Hom^\D(\tau,Q)$, for $Q\in \M^\D$. Unitality and 
associativity of the action $\kappa_Q$, for any $Q\in \M^\D$,
imply counitality and coassociativity of the left $\C$-coaction $\tau$,
respectively. 
\end{proof}

Consider a $(\C,\D)$-bicomodule $N$, over an $A$-coring $\C$ and a $B$-coring
$\D$. 
 A $\kC$-module map $g:(L,\alpha_L)\to (M,\alpha_M)$ is said to be 
{\em ${}_\kC N$-pure} provided the sequence
$$\xymatrix{ 0\ar[r] & \ker g\,\otc  N \ar[r] &  
  L\otc N \ar[rr]^{g\otc \id_N}& & M\otc N},$$
is exact (in $\M_B$).

 \begin{thm}\Label{lemma.coten.contra}{\bf Some tensor relations.} 
 Let $(N,\Nro)$ be a left $\cC$-comodule. Then:
 \begin{zlist}
\item  For any right $A$-module  $X$,
 $\rhoma \cC X\otc N \iso X\ota N.$
\item If  $(M,\roM)$ is right $\cC$-comodule for which the map 
 $$\gamma: \rhoma \cC M\to \rhoma \cC {M\ota \cC}, \; 
f\mapsto \roM\circ f - (f\ota \id_\cC)\circ \DC , $$
  is ${}_\kC N$-pure, then
 $\Rhom \cC \cC M \otc N$ is isomorphic to the cotensor product $M\ot^\C  N.$
 \end{zlist}
 \end{thm}
 \begin{proof}
(1) This is mentioned in \cite[5.1.1]{Pos}.
It is a special case of a natural isomorphism of left adjoint functors,
recalled in \ref{F-Gal}.  
Explicitly, we may
put $M=\Hom_A(\C, X)$ and $\cD=B=\End^\C(N)$ in the adjointness isomorphism to
get for $P\in \M_B$,
$$\begin{array}{rcl}
 \Hom_B(\Hom_A(\C, X)\otc N,P) &\iso &  \Hom_\kC(\Hom_A(\C, X),\Hom_B(N,P))
 \\  
      &\iso & \Hom_A(X,\Hom_B(N,P)) \\
      &\iso & \Hom_B(X\ot_A N, P).
\end{array}$$
By the Yoneda Lemma this implies the isomorphism claimed.
\smallskip

(2) Consider the commutative diagram in $\M_B$ (for $B=\End^\C(N)$).
\begin{small}
$$
\xymatrix{ 0 \ar[r] & \Rhom \cC \cC M \otc N \ar[r] 
  \ar@{-->}[d]^\vartheta & \rhoma \cC M \otc N \ar[rr]^-{\gamma\ot_\kC \id_N} 
  \ar[d]^\iso && \rhoma \cC {M\ota \cC} \otc N \ar[d]^\iso\\
  0 \ar[r] &  M \coten \cC N \ar[r]  &  M \ota N \ar[rr] &&M\ota \cC \ota N,}
$$
\end{small}
Since $\gamma$ is ${}_\kC N$-pure, and $\Rhom \cC \cC M = \ker\gamma$, 
the top row is exact. The bottom row is the defining exact sequence of the
cotensor product (see e.g.\ \cite[21.1]{BW}).
The vertical isomorphisms are obtained from part (1). Thus 
there is an isomorphism $\vartheta:\Rhom \cC \cC M \otc N\to M \coten \cC N$ 
extending the diagram commutatively.
\end{proof}

 {}From previous considerations we obtain the following result by Positselski.

\pagebreak[3]  

\begin{thm}\Label{corresp}{\bf Correspondence of categories.}
\begin{zlist}
\item For any $A$-coring $\C$, there is an adjoint pair of functors 
$$
- \otc\, \C : {\M_\kC}\to {\M^\C},\quad \Hom^\C(\C,-): {\M^\C} \to {\M_\kC}.
$$
\item For any $X\in \M_A$,
  $$ \begin{array}{rl}
X\ot_A \C \mapsto  \Hom^\C(\C,X\ot_A \C) \iso \Hom_A(\C,X), \\[+1mm]
   \Hom_A(\C,X) \mapsto  \Hom_A(\C,X)\ot_\kC \C  \iso X\ot_A \C.
\end{array}$$
Thus the functors in part (1) resrict to inverse isomorphisms 
between the {\em Kleisli subcategories} of \; $ \M^\C$ and $\M_\kC$. 
 
\item There is an equivalence
$$
\Hom^\cC(\cC,-): \M^\cC_{inj} \to \M_{[\cC,-]}^{proj},
$$
where $\M^\cC_{inj}$ denotes the full subcategory of $\M^\cC$ of objects
relative injective to the forgetful functor $\M^\cC\to \M_A$, and 
$\M_{[\cC,-]}^{proj}$ the full subcategory of $\M_{[\cC,-]}$ of
objects relative projective to the forgetful functor 
$\M_{[\cC,-]}\to \M_A$. 
\end{zlist}
\end{thm}
\begin{proof}
This is shown in \cite[Theorem in 5.3]{Pos}.
Here (1) follows by putting in \ref{thm:EM} $\cD=\cC$ and considering $\C$ as
a $(\C,\C)$-bicomodule, see \ref{HomC}. 
Claim (2) (cf. \ref{Klei}) is obtained by applying \ref{adj-mon}(2) to the
adjoint comonad-monad functor pair $(-\ota \C,\Hom_A(\C,-))$. 
Part (3) follows from \ref{Karoubi}.
Note that the equaliser in the more general situation of \ref{Karoubi} yields
here the equivalence functor $E(M,\varrho^M) =
\Hom^\cC(\cC,M)$, for any $(M,\varrho^M)\in \M^\cC_{inj}$, as stated. 
\end{proof}

Recall that an $A$-coring $\cC$ is said to be a {\em coseparable coring} if
its coproduct is a split monomorphism of $\C$-bicomodules. Equivalently,
there is an $A$-bimodule map 
$\delta: \cC\ota \cC\to A$ such that $\delta\circ\DC = \counit$
and 
\begin{equation*}
(\id_\cC\ota  \delta) \circ (\DC\ota  \id_\cC) = (\delta\ota \id_\cC)\circ
  (\id_\cC\ota  \DC). 
\end{equation*}
Such a map $\delta$ is called a {\em cointegral} (e.g.\ \cite[26.1]{BW}). 
Equivalently, coseparable corings can be described by  separable
 functors as follows.

\begin{thm}\Label{Cosep} {\bf Coseparable corings.} 
For $\cC$ the following are equivalent.
\begin{blist}
\item $\cC$ is a coseparable coring;
\item the forgetful functor $U^\cC:\M^\cC\to \M_A$ is separable;
\item the forgetful functor $U_{[\cC,-]}:\M_{[\cC,-]} \to \M_A$ is separable.
\end{blist}
If these assertions hold then, in particular, any $\cC$-comodule is
$U^\cC$-injective and any $[\cC,-]$-module is $U_{[\cC,-]}$-projective. 
\end{thm}

\begin{proof}
Equivalence (a)$\Leftrightarrow$(b) is quoted from \cite[26.1]{BW}. It can be
derived alternatively from \ref{sep.mon}(2). 
Equivalence (b)$\Leftrightarrow$(c) and the final claims follow by
\ref{adj-sep}(2). 
\end{proof}

Combining \ref{Cosep} with \ref{corresp} we obtain:

\begin{thm}\Label{cosep.cor} {\bf  Comodules and contramodules of 
coseparable corings.} 
For a coseparable coring $\cC$, the category $\M^\cC_{inj}$
coincides with $\M^\cC$ and $\M_{[\cC,-]}^{proj}$ is equal to
$\M_{[\cC,-]}$. Thus there is an 
equivalence $$\Hom^\cC(\cC,-): \M^\cC \to \M_{[\cC,-]}.$$
\end{thm}

This equivalence between comodules and contramodules for coseparable
corings plays an important role in the characterisation of categories of
Hopf (contra)modules in \ref{Hopf.sep}.

\section{Galois bicomodules} 
 
In this section we analyse, when a comodule category is equivalent to 
a contramodule
category. Any such equivalence is necessarily given by functors associated to
a bicomodule. The latter must possess additional properties. 

\begin{thm} \Label{can-iso}{\bf \CGbcm s.}
  For a $(\C,\D)$-bicomodule $(N,\Nro,\roN)$,
  the commutative diagram in \ref{contra.prod} yields 
 a canonical monad morphism (by \ref{F-Gal})
$$  
\can^N: \Hom_A(\C,-) \to \Hom^\D(N, - \otimes _A N),
  \quad   f\mapsto (f \ota \id_N) \circ\Nro.
$$  
Let $\eta$ denote the unit of the adjunction $(-\ot_\kC N,\Hom^\D(N,-))$ in
\ref{thm:EM}.  
 
The following statements are equivalent:

\begin{blist}
\item  The natural transformation $\can^N$  is an isomorphism;
\item  $\eta_{\Hom_A(\C,Q)}$ is an isomorphism, for all $Q\in
  \M_A$;
\item $\eta_M$ is an isomorphism, for all $U_\kC$-projective  
    $M\in \M_\kC$.
\end{blist}

\em If these conditions hold,
then $\Hom^\cD(N,-)$ is a  $[\C,-]$-Galois functor and we call 
$N\in {}^\C\M^\D$ a {\em \CGbcm}.
\end{thm}

If $N$ is just a left $\C$-comodule we tacitly take $\cD=B={\End^\C}(N)$ 
and call $N$ a {\em \CGlcm}. 
 
Symmetrically to the above considerations, any right adjoint functor from the
category  of left comodules of a $B$-coring $\D$ to the category of left
contramodules of an $A$-coring $\C$ is naturally isomorphic to $\Hom^\D(N,-)$,
for some $(\D,\C)$-bicomodule $N$. 
In analogy with \ref{can-iso}, also 
{\em $\Hom_{A,-}(\C,-)$-Galois $(\D,\C)$-bicomodules} and in particular    
{\em $\Hom_{A,-}(\C,-)$-Galois right $\C$-comodules} can be defined.

Studying \CGbcm s we are on the one side interested 
in there own structural properties and on the other side also in 
conditions which make the related functors fully faithful.

\begin{thm}\Label{thm:left-ff}{\bf $-\ot_\kC N$ fully faithful.}
Let $N$ be a $(\C,\D)$-bicomodule. Then
the functor $-\ot_\kC N:\M_\kC \to \M^\D$ is fully
faithful if and only if 
\begin{rlist}
\item  $N$ is a \CGbcm\  and 
\item for any $[\C,-]$-module $M$, 
the functor $\Hom^\D(N,-):\M^\D \to {\M}_A$ 
      preserves the coequaliser  
\[ \xymatrix{
\Hom_A(\C,M)\otimes_A N
\ar@<.5ex>[rr]
\ar@<-.5ex>[rr]
&&
M\otimes_A N
\ar[rr]
&&
M\ot_\kC N\ ,
} \]
defining the contratensor product (cf. \ref{contra.prod}).
\end{rlist}
\end{thm}

\begin{proof}
Since in $\M^\D$ any parallel pair of morphisms has a coequaliser,
the claim follows by (the dual version of) 
\cite[Theorem 2.6]{Pepe:Comonad&coring}.  
\end{proof}

\begin{corollary}\Label{cor:N-proj}
Let $N\in {^\C\M^\D}$ be a \CGbcm .  
If the functor $\Hom^\D(N,-):{\M}^\D \to \M_A$ preserves coequalisers,
then $-\ot_\kC N:\M_\kC \to \M^\D$ is fully faithful and $\C$ is a
projective right $A$-module. 
\end{corollary}

\begin{proof}
It follows immediately by \ref{thm:left-ff} that $-\ot_\kC
N:\M_\kC \to \M^\D$ is fully faithful. 
The left adjoint functor $-\ot_A N$  always preserves cokernels and
$\Hom^\D(N,-)$ does so by hypothesis. Thus their composite $\Hom^\D(N,-\ota
N):\M_A \to \M_A$ preserves cokernels, i.e.\ epimorphisms.
Since $\can^N$ in \ref{can-iso} is assumed to be an isomorphism, we conclude
that also the functor $\Hom_A(\C,-)$ preserves
epimorphisms, i.e.\ $\C$ is projective as a  right $A$-module. 
\end{proof}

For our investigation it is of interest to extend 
the notion of Galois comodules from \cite[4.1]{WisGal}  
to bicomodules.

\begin{thm}\Label{Galois-com}{\bf $\D$-Galois bicomodules.} \em
For any $(\C,\D)$-bicomodule $N$, the left adjoint functor    
$- \ot_\kC N:\M_\kC\to \M_B$ is a left $-\ot_B\D$-comodule functor 
(in the sense of \cite[3.3]{MW}) by the coaction 
$$-\ot_\kC \varrho^N: - \ot_\kC N \to - \ot_\kC N\ot_B\D.$$
We call $N$ a {\em \DGbcm}\ 
if $- \ot_\kC N:\M_\kC\to \M_B$ is a $-\ot_B\D$-Galois functor (in the sense
of \cite[3.3]{MW}), that is, if the comonad morphism
\begin{small}
$$
\xymatrix{
\Hom_B(N,-)\ot_\kC N \ar[rrr]^-{\id_{\Hom_B(N,-)}\ot_{[\C,-]} \varrho^N}&&&
  \Hom_B(N,-)\ot_\kC N\ot_B \D
   \ar[r]^-{\adjcu \ot_B \id_\D} & -\ot_B\D }
$$ 
\end{small}
is an isomorphism. 
\end{thm}

For a right $\D$-comodule $N$, one can put $\C=A= \End^\D(N)$. In this way we
re-obtain the usual notion of a \DGrcm\ in
\cite[4.1]{WisGal}. 

The $\D$-Galois property of a $(\D,\C)$-bicomodule is defined symmetrically  by
the Galois property of the induced functor between the category of left
$\D$-comodules and the category of left $\C$-contramodules. In the particular
case of a left $\D$-comodule $N$, it reduces to the usual notion of a
\DGlcm\ in \cite{BW} by putting $A=\C=\End^\D(N)$. 

\begin{thm}\Label{Hom-ff}{\bf $\Hom^\D(N,-)$ fully faithful.}
Let $N$ be a $(\C,\D)$-bicomodule. Then
the functor $\Hom^\D(N,-):\M^\D\to M_\kC $ is fully
faithful if and only if 
\begin{rlist}
\item  $N$ is a \DGbcm\ and 
\item the functor 
$-\ot_\kC N:\M_\kC \to \M_B$ preserves the equaliser 
\[
\xymatrix{
\Hom^\D(N,Q) \ar[rr]
&& \Hom_B(N,Q)
 \ar@<.5ex>[rr]^-{\Hom_B(N,\varrho^Q)}
\ar@<-.5ex>[rr]_-{\omega}&&
 \Hom_B(N,Q\ot_B \D),
}
\]
for any right $\D$-comodule $(Q,\varrho^Q)$, where
$\omega(f)=(f\ot_B \id_\D)\circ \varrho^N$.
\end{rlist}
\end{thm}

\begin{proof}
 This follows again by \cite[Theorem 2.6]{Pepe:Comonad&coring}.  
\end{proof}

\begin{corollary}\Label{cor:D-flat}
Let $N\in {^\C\M^\D}$ be a \DGbcm .  
If the functor $-\ot_\kC N:\M_\kC \to \M_B$ preserves 
equalisers,
then $\Hom^\D(N,-):{\M}^\D \to \M_\kC$ is fully faithful and
$\D$ is a flat left $B$-module. 
\end{corollary}

\begin{proof} (Compare with \cite[4.8]{WisGal}).
The first assertion follows immediately from \ref{Hom-ff}.
In particular, this means that $N$ is a generator in $\M^\D$. 
Moreover, there is a natural isomorphism $-\ot_B\D\iso \Hom_B(N, -)\ot_\kC N$,
where the right adjoint functor $\Hom_B(N,-)$ 
always preserves kernels and by assumption so does 
$-\ot_\kC N$. This implies that $-\ot_B \D:\M_B \to \M_B$ preserves kernels,
i.e.\ monomorphisms, hence $\D$ is flat as a left $B$-module. 
\end{proof}

Recall from \cite[19.19]{BW}
that, for a left $\C$-comodule $(N,{^N\!\varrho})$ which is 
finitely generated and projective as a left $A$-module, the left dual
${}^*N=\Hom_{A,-}(N,A)$ carries a canonical right $\C$-comodule structure, via  
$${}^* N \lra \Hom_{A,-}(N,\C) \iso {}^* N \ot_A \C, \quad 
  g \mapsto (\id_\C\ot_A g) \circ {^N\!\varrho}.$$

In what follows, $[\C,-]$-Galois and $\C$-Galois properties
of a finitely generated projective comodule are compared.
 
\begin{thm} \Label{lemma.gal.fin}{\bf $[\C,-]$-Galois comodules and
    $\C$-Galois comodules.}
 Let $N$ be a left $\cC$-comodule finitely generated and projective as a left
 $A$-module. The following assertions are equivalent.
\begin{blist}
\item $N$ is a $\Hom_{-,A}(\C,-)$-Galois left comodule;
\item $N$ is a $\C$-Galois left comodule;
\item ${}^*N$ is a $\Hom_{A,-}(\C,-)$-Galois right comodule;
\item ${}^*N$ is a $\C$-Galois right comodule.
\end{blist}
\end{thm}
 
\begin{proof} (b)$\LRa$(d) is proven in \cite[p 514]{Brz:GalCom}.  

(a)$\LRa$(d). Put $B=\End^\C(N)$ and consider the $(A,B)$-bimodule $N$ and
the $(B,A)$-bimodule ${}^*N$.
The stated equivalence follows by applying \ref{quest}(2) to the adjoint
comonad-monad pair $(-\ota\C,\Hom_A(\C,-))$ and the functor $-\ota N\iso
\Hom_{A}({}^*N,-):\M_A \to \M_B$, possessing the right adjoint $\Hom_B(N,-)$ and
the left adjoint $-\ot_B {}^* N$.  

(b)$\LRa$(c) is proven similarly to (a)$\LRa$(d).
\end{proof}

Sufficient and necessary conditions for the equivalence between a comodule and
a contramodule category are obtained by applying Beck's theorem; see
\ref{Beck}.   

\begin{thm}\Label{thm:PTT}{\bf Equivalences.}
For an $A$-coring $\C$ and a $B$-coring $\D$, the following
assertions are equivalent.

\begin{blist}
\item The categories $\M_\kC$ and $\M^\D$ are equivalent;
\item there exists a $(\C,\D)$-bicomodule $N$ with the properties:
 \begin{rlist}
 \item  $N$ is a \CGbcm ,
 \item  the functor $\Hom^\D(N,-):\M^\D \to \M_A$ reflects isomorphisms,
 \item the functor $\Hom^\D(N,-):\M^\D \to \M_A$ preserves 
       $\Hom^\D(N,-)$-contractible coequalisers. 
\end{rlist}

\item there exists a $(\C,\D)$-bicomodule $N$ with the properties:
\begin{rlist}
\item  $N$ is a \DGbcm ,
\item  the functor $-\ot_\kC N:\M_\kC \to \M_B$ reflects isomorphisms,
\item the functor $-\ot_\kC N:\M_\kC \to \M_B$ preserves 
      $-\ot_\kC N$-contractible equalisers. 
\end{rlist}
\end{blist}
\end{thm}

\begin{proof}
(a)$\LRa$(b). 
By \ref{thm:EM}, any equivalence functor $M^\D \to \M_\kC$ is naturally
isomorphic to $\Hom^\D(N,-)$, for some $(\C,\D)$-bicomodule $N$.
By Beck's theorem \ref{Beck}, $\Hom^\D(N,-):M^\D \to \M_\kC$ is an equivalence
if and only if the conditions in part (b) hold.

(a)$\LRa$(c) is shown with similar arguments.   
\end{proof}

\begin{thm}\Label{cor:str-thm}{\bf Equivalence for abelian categories.}
For a $(\C,\D)$-bicomodule $N$, the
following   are equivalent.
\begin{blist}
\item  $\Hom^\D(N,-):\M^\D \to \M_\kC$ 
is an equivalence, $\C$ is a projective right $A$-module and
  $\D$ is a flat left $B$-module;
\item $\D$ is flat as a left $B$-module and  $N$ is a \CGbcm\ and a   
  projective generator in $\M^\D$;
\item $\C$ is projective as a right $A$-module and $N$ is a \DGbcm\ and the
  functor 
  $-\ot_\kC N:\M_\kC \to \M_B$ is left exact and faithful.
\end{blist}
\end{thm}

\begin{proof}  
(a)$\Rightarrow$(b). By Theorem \ref{thm:PTT},
 $N$ is a \CGbcm . Being an equivalence, 
  $\Hom^\D(N,-):\M^\D\to\M_\kC$ is faithful. Since the forgetful functor 
  from $ \M_\kC$ to $\M_A$ (or to $\M_\Z$ or $\mathrm{Set}$) is faithful, 
 so is the composite
 $\Hom^\D(N,-): \M^\D \to \mathrm{Set}$. This proves that $N$
 is a generator in $\M^\D$. Finally, $U_\kC: \M_\kC \to
 \M_A$ is right exact by  \ref{MC} (8). Since
 $\Hom^\D(N,-):\M^\D \to \M_\kC$ is an equivalence,
 this implies that also $\Hom^\D(N,-):\M^\D \to {\M}_A$  
is right exact, by commutativity of the   diagram
in \ref{contra.prod}. By flatness of $\D$ as a left $B$-module, this 
implies projectivity of $N$ (cf. \cite[18.20]{BW}).

(b)$\Rightarrow$(a). By the hypothesis,
the functor $\Hom^\D(N,-):\M^\D \to \M_A$ preserves
coequalisers   and reflects isomorphisms.
 Thus $\Hom^\D(N,-):\M^\D \to {\M}_\kC$ 
is an equivalence by Theorem \ref{thm:PTT} and $\C$ is a projective right
$A$-module by Corollary \ref{cor:N-proj}. 

(a)$\Rightarrow$(c). If $\Hom^\cD(N,-):\M^\cD \to \M_{[\cC,-]}$ is an
  equivalence, then so is its left adjoint $-\otimes_{[\cC,-]} N$. Thus $N$ is
  a \DGbcm\ by \ref{thm:PTT}. The functor $-\otimes_{[\cC,-]} N:
  \M_{[\cC,-]} \to \M_B$ is equal to the composite of the equivalence
  $-\otimes_{[\cC,-]} N: \M_{[\cC,-]} \to \M^\cD$ and the forgetful functor
  $\M^\cD\to \M_B$. The forgetful functor is faithful and also left exact by
  the flatness of the left $B$-module $ \cD$. Thus the functor 
$-\otimes_{[\cC,-]} N: \M_{[\cC,-]}
  \to \M_B$ is also faithful and left exact. 

(c)$\Rightarrow$(a). Since $\cC$ is a projective right $A$-module,
  $\M_{[\cC,-]}$ is 
  abelian. Hence faithfulness of $-\otimes_{[\cC,-]} N: \M_{[\cC,-]} \to \M_B$
  implies that it reflects isomorphisms. Since it also preserves equalisers by
  assumption, it follows by Theorem \ref{thm:PTT} that 
  $-\ot_{[\cC,-]} N:\M_{[\cC,-]} \to
  \M^\cD$ is an equivalence, with inverse $\Hom^\cD(N,-)$. The left $B$-module
  $\cD$ is flat by Corollary \ref{cor:D-flat}.
\end{proof}

In the rest of the section we study the particular case of a trivial
$B$-coring $\D=B$. That is, the situation when the category of contramodules of
a coring $\C$ is equivalent to that of modules over a ring $B$.

\begin{lemma}\Label{lem:T-B}\cite[Proposition 2.5]{ElKaGomTor:Comat}
Let $N$ be an $(A,B)$-bimodule which is
finitely generated and projective as an $A$-module. Consider the comatrix
coring $\C:=N\otimes_B {}^* N$ and denote by $T$ the ring of endomorphisms of
$N$ as a left $\C$-comodule.
Then $N\iso N\otimes_B T$ via the right $T$-action on $N$. 
\end{lemma}

The next result may be seen as a counterpart to the Galois comodule structure
theorem \cite[18.27]{BW}.

\begin{theorem} \Label{thm:(contra)Gal}
 Let $N\in {^\C\M}$ be a \CGcm\ over an $A$-coring $\C$, put 
 $T=\End^\C(N)$ and $B\subset T$ be a subring.
Assume that $N$ is a projective generator of right
$B$-modules. Then the following hold.
\begin{zlist}
\item   $-\ot_\kC N:\M_\kC\to \M_B$ is an equivalence.
\item  $\C$ is a projective right $A$-module.
\item  $N$ is a finitely generated and projective left $A$-module. 
\item  $\C$ is isomorphic to the comatrix $A$-coring $N \ot_B {}^*N$. 
\item  $B$ is isomorphic to $T$.
\item  If, in addition, $\C$ is a generator of right $A$-modules, then $N$
  is a faithfully flat left $A$-module.  
\end{zlist}
\end{theorem}

\begin{proof}
Assertions  (1) and (2)  are immediate by  \ref{cor:str-thm}. 

 (3)   Since $-\ot_\kC N$ is an equivalence, it has a left
adjoint $\Hom_B(N,-):\M_B \to \M_\kC$. The free
functor $\Hom_A(\C,-)$ has a left adjoint $-\ot_\kC \C:\M_\kC
\to \M_A$  by  \ref{contra.prod}. 
Hence also the composite functor, 
that is naturally isomorphic to $ - \otimes_A N:\M_A \to
\M_B$ by \ref{lemma.coten.contra}(1),
has a left adjoint. This
proves that $N$ is a finitely generated and projective left $A$-module.  

 (4) By part (3), 
\[
\Hom_B(N, - \otimes_A N)\iso 
\Hom_B(N,\Hom_A({}^* N,-)) \iso 
\Hom_A(N\otimes_B {}^* N,-).
\]
Composing this natural isomorphism with 
the canonical monad morphism $\can^N$ (at $\D=B$), it yields a
monad isomorphism $\Hom_A(\C,-)\iso \Hom_A(N\otimes_B {}^*
N,-)$. By Yoneda's Lemma this proves $\C\iso N \otimes_B {}^* N$.
\smallskip

 (5)  The composite of the forgetful functor
$\Hom_T(T,-): \M_T \to \M_B$ and
$\Hom_B(N,-): \M_B \to \M_A$ is naturally
isomorphic to 
\[
\Hom_B(N,\Hom_T(T,-))\iso \Hom_T(N\otimes_B T,-)\iso  \Hom_T(N,-),
\]
where the last isomorphism follows by part (4) and Lemma \ref{lem:T-B}.
The forgetful functor $\M_T \to \M_B$ reflects
isomorphisms. Since $N$ is a generator in $\M_B$ by  
assumption,  the (fully faithful) functor $\Hom_B(N,-):
\M_B \to \M_A$ reflects isomorphisms too. Hence also the composite
$\Hom_T(N,-):\M_T \to \M_A$ reflects
isomorphisms. The forgetful functor $\M_T \to \M_B$ has a
right adjoint (the coinduction functor $\Hom_B(T,-)$) hence it
preserves coequalisers. Since $N$ is a projective right $B$-module by
assumption,  $\Hom_B(N,-): \M_B \to \M_A$ preserves coequalisers too.
Hence also the composite 
$\Hom_T(N,-):\M_T \to {\M}_A$ 
preserves coequalisers. The equivalence functor 
$-\ot_\kC N: {\M}_\kC \to \M_B$ factorises through 
$-\ot_\kC N: {\M}_\kC \to \M_T$ 
and the forgetful functor $\M_T \to
\M_B$. Thus the forgetful functor is full (and obviously
faithful). This implies that $-\ot_\kC N: \M_\kC \to \M_T$
is fully faithful, hence the corresponding canonical monad morphism
$$
\Hom_A(\C,-) \to \Hom_T(N, - \otimes _A N),
  \quad f \mapsto [\  n\mapsto (f\ota \id_N) \circ{}^N\varrho(n)] ,
$$
is a natural isomorphism by Theorem \ref{thm:left-ff}. So we conclude by
Theorem \ref{thm:PTT} that $-\ot_\kC N: \M_\kC \to \M_T$
is an equivalence and so is the forgetful functor $\M_T \to
\M_B$. This proves the isomorphism of algebras $T\iso B$.
\smallskip

 (6)  $N$ is a flat left $A$-module by part (3). Hence it
suffices to show that, under the assumptions made, $-\otimes_A N: 
\M_A \to \M_B$ is a faithful functor, so it reflects both
monomorphisms and epimorphisms. Recall that, by 
\ref{lemma.coten.contra}(1),
$-\otimes_A N: \M_A \to \M_B$ is naturally isomorphic to
the composite of the free functor $\Hom_A(\C,-): \M_A \to
\M_\kC$ and the equivalence $-\ot_\kC N:\M_\kC \to 
\M_B$. By assumption, $\Hom_A(\C,-): \M_A \to 
\M_A$ is faithful. Then also $\Hom_A(\C,-): \M_A \to 
\M_\kC$ is faithful, what completes the proof.
\end{proof}

Note that $\C$ is a generator of right $A$-modules as in Theorem
\ref{thm:(contra)Gal} (6) in various situations. For example, whenever the
counit of 
$\C$ is an epimorphism (e.g.\ because there exists a grouplike element in $\C$
or $\C$ is faithfully flat as a left or right $A$-module).

 \section{Contramodules and entwining structures}

As recalled in \ref{lift.e}, lifting of a monad
$\bF$ on a category $\A$ to a monad on the category 
$\A^\nG$ for a comonad $\bG$, 
or lifting of a comonad $\bG$ to a comonad on the
category $\A_\nF$ for a monad $\bF$, are both equivalent to the
existence of a mixed distributive law (entwining) between $\bF$ and
$\bG$. Combining this general fact with properties of module
categories, we obtain a description of entwinings between 
$A$-rings and $A$-corings ($A$ is an associative ring with unit). 
Recall that a {\em (left) entwining map} between an $A$-ring $B$ and
an $A$-coring $\C$ is an $\AA $-bimodule morphism
$\psi:B \otimes_A \C \to \C \otimes_A B$ which respects (co)multiplications 
and (co)units. Similarly, (right) entwining maps 
$\lambda: \C\otimes_A B\to B\otimes_A \C$ 
are defined (e.g.\  \cite[Chapter~5]{BW}).

\begin{thm} \Label{thm:entwining}{\bf Entwining maps.}
For  all $A$-rings  $B$  and $A$-corings $\C$, the following
assertions are equivalent.
\begin{blist}
\item There is an entwining map $\psi:B \otimes_A \C \to \C \otimes_A B$;
\item the monad $B\otimes_A -$ on ${}_A\M$ has a lifting
      to a monad on   ${}^\C \M$;
\item the comonad $\C\otimes_A -$ on ${}_A\M$ has a
      lifting to a comonad on ${}_B {\M}$;
\item the monad $\Hom_A(\C,-)$ on $\M_A$ has a
      lifting to a monad on ${\M}_B$; 
\item the comonad $\Hom_A(B,-)$ on $\M_A$  has a
      lifting to a comonad on  $\M_\kC$.
\end{blist}
\end{thm}

\begin{proof}
(a)$\Leftrightarrow$(b)  and (a)$\Leftrightarrow$(c). 
An entwining map $\psi$ determines a mixed distributive law $\Psi:=
\psi \otimes_A -:B \otimes_A \C \otimes_A - \to \C \otimes_A B \otimes_A
-$. Conversely, if $\Psi: B \otimes_A \C \otimes_A - \to \C \otimes_A B
\otimes_A -$ is a mixed distributive law, then $\psi:=\Psi_A$ is an entwining
map. 
\smallskip

(a)$\Leftrightarrow$(d)  and (a)$\Leftrightarrow$(e). 
An entwining map $\psi$ determines a mixed distributive law ${\widetilde \Psi}$:
$$\begin{array}{rcl}
 \Hom_A(\C, \Hom_A(B,-))\iso\quad    && \quad\iso \Hom_A(B, \Hom_A(\C,-)) \\
  \Hom_A(\C \otimes_A B,-)   &\stackrel{\Hom_A(\psi,-)}\lra & 
  \Hom_A(B \otimes_A \C,-).  
\end{array}$$

On the other hand, by the Yoneda Lemma, any mixed distributive law 
${\widetilde \Psi}:
\Hom_A(\C \otimes_A B,-) \to \Hom_A(B \otimes_A \C,-)$ is of
this form.
\end{proof}

Under the equivalent conditions of \ref{thm:entwining}, $\C\otimes_A B$
is a $B$-coring, cf. \cite[32.6]{BW}. 
Its contramodules can be described as follows. 

\begin{thm}\Label{Mod-str}{\bf $[\C\otimes_A B,-]$-modules.}
Let $B$ be an $A$-ring and 
 $\C$ an $A$-coring with an entwining map  $\psi:B
\otimes_A \C \to \C \otimes_A B$. Then the following
structures on a right $B$-module $M$ are equivalent.
\begin{blist}
\item A module  structure map $\varrho_M :\Hom_B (\C\otimes_A B,M) \to M $;
\item a   $B$-linear module structure map
 $\alpha_M :\Hom_A (\C,M) \to M $
  (where $fb=\sum f(-^\psi) b_\psi$, for $f\in \Hom(\C,M)$, $b\in B$, 
  hence $B$-linearity means
  $\alpha_M(f)b = \sum \alpha_M\big(f(-^\psi) b_\psi\big)$, 
  with notation  
  $\psi(b\otimes_A c)= \sum c^\psi \otimes_A b_\psi$);
\item  a module structure for the monad $\Hom_A(\C,-)$ on
  $\M_B$;
\item a comodule structure for the comonad $\Hom_A(B,-)$ on $\M_\kC$. 
\end{blist}
\end{thm}

\begin{proof}
(a)$\Leftrightarrow$(b). The isomorphism
$\Hom_B(\C\otimes_A B,M)\iso \Hom_A(\C,M)$ of right $A$-modules
induces an isomorphism 
\[
\xi:\Hom_A(\Hom_A(\C,M),M) \to
\Hom_A(\Hom_B(\C\otimes_A B,M),M).
\] 
As easily checked, $\xi (\alpha_M)$ belongs to
$\Hom_B(\Hom_B(\C\otimes_A B,M),M)$ if and only if  $\alpha_M$
satisfies the $B$-linearity condition in (b). Associativity and unitality of a
$[\C\otimes_A B,-]$-action $\xi (\alpha_M)$ are equivalent to analogous
properties of the $[\C,-]$-action $\alpha_M$.

Equivalences (b)$\Leftrightarrow$(c)  and 
(b)$\Leftrightarrow$(d)  follow by \ref{f.entw} (cf. \cite[5.7]{WisAlg}).
\end{proof}

In light of \ref{thm:entwining}, the following describes a special case of
\ref{dist-adj} and \ref{dist-mod}.

 \begin{thm}\Label{dist-AC}{\bf Distributive laws for rings and corings.} 
  Let $B$ be an $A$-ring and $\C$  an $A$-coring over any ring $A$.
\begin{zlist}
\item $\lambda:\C\ot_A B\to B\ot_A\C$ is an entwining map if and only if 
$$ 
\tilde \lambda: \Hom_A(B,-)\ot_A \C\to \Hom_A(B, - \ot_A \C),\;
f \ot c \mapsto (f \ot_A\id_\C) \circ \lambda(c\ot -),
$$
is a comonad distributive law. 
In this case, $\Hom_A(B,-)\ot_A \C$ is a
comonad on $\M_A$ and the category of its comodules is isomorphic to the
category of $-\ot_A \lambda$-bimodules (i.e.\ usual entwined modules),
cf.\ \ref{f.entw}. 
 
\item $\psi:B\ot_A \C\to \C\ot_AB$ is an entwining map if and only if 
$$
\tilde\psi: \Hom_A(\C,-)\ot_A B \to \Hom_A(\C,-\ot_AB),\; 
g\ot a \mapsto (g\ot_A \id_B)\circ \psi(a\ot -) ,
$$
is a monad distributive law. 
In this case, $\Hom_A(\C,-\ot_AB)$ is a monad on $\M_A$ and the category of
its modules is isomorphic to the category of $\Hom_A(\psi,-)$-bimodules
(cf.\ \ref{f.entw}). 
\end{zlist}
\end{thm}

Note that for a commutative ring $R$, any $R$-algebra $A$ and $R$-coalgebra
$C$ are entwined by the twist 
maps $C\ot_R A \to A\ot_R C$ and $A\ot_R C\to C\ot_R A$. Applying
\ref{dist-AC} to these particular entwinings, we conclude that the
canonical natural transformations
\begin{eqnarray*}
&\Hom_R(A,-) \ot_R C \to \Hom_R(A,-\ot_R C), \quad &f\ot c\mapsto  f(-)\ot c, 
\quad \mathrm{and}\\
&\Hom_R(C,-) \ot_R A \to \Hom_R(C,-\ot_R A), \quad &g\ot a\mapsto  g(-)\ot a,
\end{eqnarray*}
yield a comonad distributive law and a monad distributive law, respectively.

\section{Bialgebras and bimodules}

There are many equivalent characterisations of bialgebras and Hopf algebras.
A bialgebra over a commutative ring $R$ can be seen as an $R$-module that is
both an algebra and a coalgebra entwined in a certain way. In category theory
terms, bialgebra is defined as an $R$-module such that the tensor functor
$-\ot_R B$ is a bimonad on $\M_R$. 
Associated to a bialgebra $B$, there is a category of Hopf modules,
whose objects are $B$-modules with a compatible $B$-comodule structure. A Hopf
algebra can be characterised as a bialgebra $B$ such that the functor 
$-\ot_R B$ is an equivalence between the categories of $R$-modules and Hopf
$B$-modules. In this section we supplement this description of bialgebras and 
Hopf algebras by  the equivalent description in terms of properties of the
Hom-functor $[B,-]$, and hence in terms of contramodules. 

Throughout, $R$ is a commutative ring. 
The unit element of a (bi)algebra $B$ is denoted by $1_B$. 
For the coproduct $\cop$ of a
bialgebra $B$, if applied to an element $b\in B$, we use Sweedler's index
notation $\cop(b)= b \sw1 \ot b \sw2$, where implicit summation is
understood. 

\begin{thm}\Label{Bialg} {\bf Bialgebras.} \em 
Let $B$ be an $R$-module which is both 

an $R$-algebra $\product:B\ot_RB\to B$, $\unit:R\to B$, and

  an $R$-coalgebra  $\cop : B\to B\ot_RB$,
 $\counit: B\to R$. \\
Based on the canonical twist $\tw:B\ot_RB\to B\ot_RB$, we obtain the following 
$R$-module maps
$$ \begin{array}{rl}
\psi_r = (\id_B\ot_R \product )\circ  (\tw\ot_R \id_B)\circ (\id_B\ot_R \cop):&
   B\ot_RB\to B\ot_RB, 
\\
\psi_l = (\product \ot_R \id_B)\circ  (\id_B \ot_R\tw)\circ (\cop\ot_R \id_B):&
   B\ot_RB\to B\ot_RB.
\end{array} $$
Evaluated on elements, $\psi_r(a\ot b)=b\sw1 \ot a b\sw2
$ and $\psi_l(a\ot b) = a\sw1 b \ot a \sw2$.

To make $B$ a {\em bialgebra}, $\product$ and $\unit$ must be coalgebra maps 
(equivalently, $\cop$ and $\counit$ are to be algebra maps) with respect to
the obvious product and coproduct on $B\ot_RB$ (induced by $\tw$).  
The compatibility between multiplication and comultiplication can be 
expressed by commutativity of the diagram
 $$\xymatrix{ B\ot_R B \ar[r]^-{\product} \ar[d]_{\cop\ot_R \id_B} &B
  \ar[r]^-{\cop} & B\ot_R B \\ 
 B\ot_R B\ot_R B\ar[rr]^{\id_B\ot_R \psi_r} & & B\ot_R B\ot_R B
 \ar[u]_-{\product \ot_R \id_B}.}
$$

For a bialgebra $B$, both maps $\psi_r$ and $\psi_l$ are (right, respectively,
left) entwining maps between the algebra $B$ and the coalgebra $B$. Going to
the functor level it turns out that $\psi_r$ yields a mixed distributive law
for the monads and comonads $-\ot_R B$, while $\psi_l$ is related to the
endofunctors $B\ot_R-$.   
\end{thm}

Given an $R$-bialgebra $B$, it will sometimes help to write $\uB$ when we
focus on the algebra structure and $\oB$ when focussing on the coalgebra part.
 From \ref{thm:entwining} we know:

\begin{thm} \label{thm:entwining.bi}{\bf Entwining maps for bialgebras.}
Consider an $R$-module $B$ which is an $R$-algebra and an $R$-coalgebra. The
following assertions are equivalent.
\begin{blist}
\item There is an entwining map $\psi:\uB \otimes_R \oB \to \oB \otimes_R \uB$;
\item the monad $\uB\otimes_R -$ on $\M_R$ has a lifting
      to a monad on   ${}^\oB \M$;
\item the comonad $\oB\otimes_R -$ on $\M_R$ has a
      lifting to a comonad on ${}_\uB {\M}$;
\item the monad $\Hom_R(\oB,-)$ on $\M_R$ has a
      lifting to a monad on ${\M}_\uB$; 
\item the comonad $\Hom_R(\uB,-)$ on $\M_R$  has a
      lifting to a comonad on  $\M_\kB$.
\end{blist}
\end{thm}
A symmetric form of Theorem \ref{thm:entwining.bi} can be obtained by
interchanging left and right (co)module structures.

\begin{thm}\Label{cor.entw}{\bf Entwinings and corings.} \em
Let $B$ be  an $R$-bialgebra with an entwining map  
$\psi:  \oB\otimes_R \uB \to \uB \otimes_R \oB$. Then
$B\ot_RB$ is a $\uB$-coring with structure maps 
$$
\id_B \ot_R \cop:  B\ot_R B\to ( B\ot_R B)\ot_B ( B\ot_R B),\quad
\id_B \ot_R \counit: B\ot_RB\to B, 
$$
and $\uB$-actions $d\cdot (a\ot b)\cdot c= da \psi(b\ot c)$. 

Symmetrically, an entwining map $\psi:\uB \ot_R \oB \to\oB \ot_R \uB$
determines a $\uB$-coring $B \ot_R B$, with structure maps
$$
\cop \ot_R \id_B: B\ot_R B\to ( B\ot_R B)\ot_B ( B\ot_R B), \quad
\counit \ot_R \id_B: B\ot_RB\to B,
$$
and $\uB$-actions $d\cdot(a\ot b)\cdot c=\psi(d\ot a)bc$.

In particular, the entwining maps $\psi_r:\oB\ot_R \uB\to \uB \ot_R \oB$ and
$\psi_l:\uB\ot_R \oB \to \oB \ot_R \uB$ in \ref{Bialg}
determine $\uB$-corings $B \ot_R B$, denoted by $B\ot^r_R B$ and 
$B\ot_R^lB$, respectively. 
\end{thm}

\begin{thm}\Label{Hopf.M}{\bf $B$-Hopf modules.}  
Let $B$ be an $R$-bialgebra and consider the $\uB$-coring $B \ot^r_R B$ in
\ref{cor.entw}. 
The following structures on a right $\uB$-module $M$ are equivalent: 
\begin{blist}
\item A right $B\ot_R^r B$-comodule structure map $\varrho^M: M \to M \ot_B
  (B \ot^r_R B)$;
\item a right $\uB$-linear $\oB$-comodule structure map $\alpha^M:M \to M \ot_R
  B$, (where $\uB$-linearity means commutativity of the diagram  
$$
\xymatrix{ M\ot_R B \ar[r]^-{\alpha_M} \ar[d]_-{\alpha^M\ot_R \id_B} &  
     M \ar[r]^-{\alpha^M} & M\ot_R B \\
   M\ot_R B \ot_R B  \ar[rr]^-{\id_M \ot_R \psi_r} &&  
    M\ot_R B \ot_R B \ar[u]_-{\alpha_M\ot_R \id_B } ,}
$$
where $\alpha_M: M \ot_R B \to M$ denotes the $\uB$-action on $M$);
\item a comodule structure for the comonad $- \ot_R B$ on $\M_\uB$;
\item a module structure for the monad $-\ot_R B$ on $\M^\oB$. 
\end{blist}
\end{thm}

A right $\uB$-module $M$ with these equivalent properties is called a 
{\em $B$-Hopf module}. Morphisms of $B$-Hopf
modules are $B\ot_R^r B$-comodule maps. Equivalently, they are $\uB$-module
as well as $\oB$-comodule maps. The category of right $B$-Hopf modules is 
denoted by $\M^\oB_\uB$. 
By the above considerations, it is isomorphic to $\M^{B\ot_R^r B}$.

Based on $\psi_l$, left $B$-Hopf modules are defined in a symmetric way.
Note that a bialgebra $B$ is both a left and a right $B$-Hopf module.  
\smallskip 

 {}From \ref{Mod-str} we obtain:

\begin{thm}\Label{Mod-str.B}{\bf $[B,-]$-Hopf modules.}  
Let $B$ be an $R$-bialgebra and consider the $\uB$-coring $B\ot^l_R B$ in
\ref{cor.entw}. Then the following structures on a right $\uB$-module $M$ are
equivalent.  
\begin{blist}
\item  
 A $[B\ot_R^lB,-]$-module structure map $\varrho_M :\Hom_B (B\otimes_R^l
 B,M) \to M $; 
\item a $\uB$-linear $[\oB,-]$-module structure map $\alpha_M :\Hom_R (B,M) \to
  M$ \\
  (i.e.\ $\alpha_M(f)b = \sum \alpha_M\big(f(b\sw1 -) b\sw2\big)$
 for $f\in \Hom_R(B,M)$, $b\in B$);
\item  a module structure for the monad $\Hom_R(B,-)$ on $\M_B$;
\item a comodule structure for the comonad $\Hom_R(B,-)$ on $\M_\kB$. 
\end{blist}
\end{thm}
A right $\uB$-module $M$ with these equivalent properties is called a
{\em $[B,-]$-Hopf module} or {\em right Hopf contramodule for $B$}. Morphisms
of $[B,-]$-Hopf modules are $B\ot_R^lB$-contra\-module maps. Equivalently,
they are $\uB$-module 
as well as $\oB$-contramodule maps. The category of $[B,-]$-Hopf modules is
denoted by $\M^{[\uB,-]}_{[\oB,-]}$. By the above considerations, it is
isomorphic to $\M_{[B\ot_R^l B,-]}$. 

Based on $\psi_r$, left Hopf contramodules modules for $B$ are defined in a
symmetric way.

Applying \ref{dist-AC}, the following alternative description of Hopf modules
is obtained. 

\begin{thm}\Label{dist-Bi}{\bf Distributive laws for bialgebras.}  
Let $B$ be an $R$-bialgebra. Then:
\begin{zlist}
\item The entwining $\psi_r$ in \ref{Bialg} induces 
a comonad distributive law 
$$\Hom_R(\uB,-)\otimes_R \oB \to   \Hom_R(\uB, -\otimes_R \oB), \quad
                   f \ot b \mapsto \sum f((-)\sw1) \ot b(-)\sw2 .$$
Hence $\Hom_R(\uB,-)\otimes_R \oB$ is a comonad on $\M_R$. 
The category of its comodules is isomorphic to the category of $B$-Hopf
modules.
\item The entwining $\psi_l$ in \ref{Bialg} induces
a monad distributive law 
$$
\Hom_R(\oB,-)\otimes_R \uB \to   \Hom_R(\oB, -\otimes_R \uB), \quad
                   f \ot b \mapsto \sum f(b\sw1-) \ot b\sw2 .
$$
Hence $\Hom_R(\oB, -\otimes_R \uB)$ is a monad on $\M_R$. 
The category of its modules is isomorphic to the category of $[B,-]$-Hopf
modules. 
\end{zlist}
\end{thm}

\begin{thm}\Label{Hopf}{\bf Hopf algebras. }\em
An $R$-bialgebra $(H,\product,\unit,\cop,\counit)$ is said 
to be a {\em Hopf algebra} if there is an
  $R$-module map $S:H \to H$, called the {\em antipode}, such that 
$$
\product \circ (\id_H \ot_R S)\circ \cop= \unit \circ \counit 
= \product \circ (S \ot_R \id_H) \circ \cop .
$$
If the antipode exists, then it is unique and it is an
anti-algebra and anti-coalgebra map.

For an $R$-Hopf algebra $H$, the $H$-corings $H \ot_R^r H$ and 
$H \ot_R^l H$ in \ref{cor.entw} are isomorphic via the mutually inverse maps 
\begin{eqnarray*}
&H \ot_R^r H \to H \ot_R^l H,\quad &a\otimes b \mapsto 
   \sum a\sw1 S(b\sw2) \otimes a\sw2 S(b\sw1) b\sw{3} ,\\ 
&H \ot_R^l H \to H \ot_R^r H,\quad &a\otimes b \mapsto 
  \sum a\sw1 S(a\sw3) b\sw1\otimes S(a\sw2) b\sw2. 
\end{eqnarray*}
\end{thm}

\begin{thm}\Label{Hopf.sep}{\bf Hopf algebras and coseparability.}
Let $H$ be an $R$-Hopf algebra. 
\begin{zlist}
\item The $H$-coring $H\ot_R^r H$ is coseparable.
\item  The following functor is an equivalence:
$$\Hom^{H\ot_R^r H}(H\ot_R^r H,-):
\M^{H\ot_R^r H} \to \M_{[H\ot_R^r H,-]}.$$ 
\item The category of $H$-Hopf modules (in \ref{Hopf.M}) and
the category of $[H,-]$-Hopf modules (in \ref{Mod-str.B}) are equivalent. 
\end{zlist}
\end{thm}
\begin{proof}
(1) A cointegral is given by 
$$
(H\ot_R ^r H)\ot_H (H\ot_R^r H) \to H, \quad 
(a\ot b)\ot_H (1_H\ot c) \mapsto aS(b)c.
$$

(2) In view of (1), this is a special case of \ref{cosep.cor}.

(3) The category of $H$-Hopf modules is isomorphic to $\M^{H\ot_R^r H}$ and
  the category of $[H,-]$-Hopf modules is isomorphic to $\M_{H\ot_R^l H}$. So
  the claim follows by coring isomorphism $H\ot_R^r H\iso H\ot_R^l H$ in
  \ref{Hopf} and part (2).  
\end{proof}

The final aim of this section is to characterise Hopf algebras via their
induced (co)monads. The following notions were introduced in \cite{WisAlg} and
\cite{MW}. Note that these terms have different meanings in  
Moerdijk \cite{Moer} and Brugui\`eres-Virelizier \cite{BruVir}. 

\begin{thm}\Label{Bimonad} {\bf Bimonads and Hopf monads. }\em
A {\em bimonad} on a category $\A$ is a functor 
  $F: {\A} \to {\A}$ with
a monad structure ${\underline F}= (F,m,i)$ and a comonad structure
${\overline F}= (F,d,e)$ subject to the compatibility conditions 
\begin{rlist}
\item $e$ is a monad morphism ${\underline F} \to \id_\A$;
\item $i$ is a comonad morphism $\id_\A \to {\overline F}$; 
\item there is a mixed distributive law $\Psi :{\underline F} {\overline F}
  \to {\overline F} {\underline F}$, satisfying 
$$d  \circ m =Fm\circ \Psi F \circ F d .$$
\end{rlist}

A bimonad $(F,m,i,d,e)$ is called a {\em Hopf monad} if there exists a
natural 
transformation $S: F \to F$, called the {\em antipode}, such that 
$$
m \circ SF \circ d = i \circ e = m \circ FS \circ d.
$$

A class of examples of bimonads is provided by the following construction in
\cite[Proposition 6.3]{MW} (see also \cite{ElKa}). Let $F$
be a functor $\A\to \A$ allowing a monad structure
${\underline F}= (F,m,i)$ as well as a comonad structure ${\overline F}=
(F,d,e)$. Consider a {\em double entwining} $\tau$, i.e.\ a natural
transformation $FF \to FF$, which is an entwining both in the sense 
${\underline F} {\overline F} \to {\overline F} {\underline F}$ and also $
{\overline F} {\underline F} \to {\underline F} {\overline F}$. The functor
$F$ is called a {\em $\tau$-bimonad} provided that the above conditions (i) and
(ii) hold and in addition 
$$
 mF\circ FFm\circ F \tau F \circ d FF\circ Fd  = d  \circ m.
$$
By \cite[Proposition 6.3]{MW}, a $\tau$-bimonad $F$ is a bimonad with respect
to the mixed distributive law $\Psi:= mF \circ F \tau\circ d  F$.

A $\tau$-bimonad with an antipode is called a {\em $\tau$-Hopf monad}.
\smallskip

As described in \cite{MW}, if a $\tau$-bimonad
$F$ has a left or right adjoint $G$, then the mates under 
the adjunction of the structure maps of the monad and comonad $F$, equip $G$
with a comonad and a monad structure, respectively. Moreover, the mate 
${\bar\tau}$ of $\tau$ under the adjunction is a double entwining for $G$,
and $G$ is a ${\bar \tau}$-bimonad. If $F$ is a $\tau$-Hopf monad, then 
$G$ is a $\bar\tau$-Hopf monad.
\end{thm}

\begin{thm}\Label{bimonad.B}{\bf The bimonad $-\ot_R B$.} \em
For an $R$-bialgebra $(B,\product ,\unit ,\cop ,\counit)$, the functor
$- \ot_R B : \M_R \to \M_R$ is a
$\tw$-bimonad, hence a bimonad with respect to the mixed distributive law
$$
(- \ot_R \id_B \ot_R \product)\circ (-\ot_R \tw \ot_R \id_B)\circ
(-\ot_R \id_B \ot_R \cop )= -\ot_R \psi_r.
$$
By duality, $\Hom_R(B,-)$ is a ${\overline \tw}$-bimonad, with coproduct
$[\product,-]$ and counit $[\unit,-]$ in \ref{mon-com}, product $[\cop,-]$ and  
unit $[\counit,-]$ in \ref{coring}, where 
$${\overline \tw}: \Hom_R(B,\Hom_R(B,-)) \to \Hom_R(B,\Hom_R(B,-))$$
 is given by switching the arguments. Thus $\Hom_R(B,-):\M_R \to \M_R$
 is a bimonad with respect to the mixed distributive law
 $\Hom_R(\psi_l,-)$:
$$\begin{array}{rcl}
\Hom_R(B,\Hom_R(B,-))\iso \quad & & \quad \iso  \Hom_R(B,\Hom_R(B,-)) \\
 \Hom_R(B\ot_R B,-) &\stackrel{\Hom_R(\psi_l,-)}\lra & \Hom_R (B\ot_R B,-).\\ 
\end{array}$$
\end{thm}

A motivating example of a ($\tw$-)Hopf monad in \cite{MW} is the functor
$-\ot_R H:\M_R \to \M_R$, induced by a Hopf algebra $H$. 
 
Summarising the preceding observations we obtain the following.
 
\begin{thm}\Label{Hopf.char}{\bf Characterisations of Hopf algebras.}
For an $R$-bialgebra $(H,\product ,\unit ,\cop,\counit)$, 
the following assertions are equivalent.
\begin{blist}
\item $H$ is a Hopf algebra;
\item the map               
 $\gamma: H\ot_R H \stackrel{\cop\ot_R \id_H}\lra H\ot_RH\ot_RH 
      \stackrel{\id\ot_R \product}\lra H\ot_R H$ is an isomorphism; 
\item $H$ is an $H\ot_R^r H$-Galois right (equivalently, left) comodule;  
\item $H$ is an $H\ot_R^l H$-Galois right (equivalently, left) comodule;  
\item $- \ot_R H$ is a $\tw$-Hopf monad on $\M_R$; 
\item for the ${\overline \tw}$-bimonad $[H,-]=\Hom_R(H,-)$, 
      the natural transformation 
$$[\gamma,-]: [H,[H,-]] 
   \stackrel{[H,[\product,-]]}{\lra} [H,[H,[H,-]]] 
    \stackrel{[\cop,[H,-]]}{\lra} [H,[H,-]]$$    
is an isomorphism; 
\item $\Hom_R(H,-)$ is a ${\overline \tw}$-Hopf monad on $\M_R$;
\item $- \ot_R H: \M_R \to \M^\oH_\uH$ is an equivalence;       
\item $\Hom_R(H,-): \M_R \to M^{[\uH,-]}_{[\oH,-]}$ is an equivalence;
\item $H$ is a $\Hom_{-,H}(H\ot_R^r H,-)$-Galois left comodule
(equivalently, a\\
$\Hom_{H,-}(H\ot_R^r H,-)$-Galois right comodule); 
\item $H$ is a $\Hom_{-,H}(H\ot_R^l H,-)$-Galois left comodule
(equivalently, a\\
$\Hom_{H,-}(H\ot_R^l H,-)$-Galois right comodule).   
\end{blist}
\end{thm}

\begin{proof}
(a)-(d) and (h) are standard equivalent characterisations of Hopf
  algebras, see e.g.\ \cite[15.2 and 15.5]{BW}. 

(a)$\LRa$(e)$\LRa$(f)$\LRa$(g) is proven in \cite{MW}.

(c)$\LRa$(j)  and (d)$\LRa$(k) follow by \ref{lemma.gal.fin}. 

(i)$\Ra$(j) follows by \ref{thm:PTT} (a)$\Ra$(c)(i).

(h)$\Ra$(i) There is a sequence of equivalences,
$$
\M^\oH_\uH \iso \M^{H\ot^r_R H} \equi \M_{[H\ot_R^r H,-]} \iso
\M_{[H\ot_R^l H,-]} \iso \M^{[\uH,-]}_{[\oH,-]},
$$
cf. \ref{Hopf.M}, \ref{Hopf.sep}, \ref{Hopf} and \ref{Mod-str.B} (note that
(h)$\Ra$(a)). Combining this composite with the equivalence in part (h), we
obtain an equivalence functor 
$$
\Hom^{H\ot_R^l H}(H \ot_R^l H,-\ot_R H)\,:\,\M_R \to \M^{[\uH,-]}_{[\oH,-]}. 
$$
We claim that the functor in part
(i) is naturally isomorphic to this equivalence, hence it is an equivalence,
too. 

The equivalence in part (h) gives rise to an $R$-module isomorphism
$$
\Hom^{H\ot_R^r H}(H\ot_R^r H,M\ot_R H)\to \Hom_R(H,M), \ 
\varPsi \mapsto (\id_M\ot_R \counit)\circ \varPsi(-\ot_R \unit),
$$
for any $R$-module $M$, that is natural in $M$. Using the coring isomorphism
$H\ot_R^r H \iso H \ot_R^l H$ in \ref{Hopf}, we can transfer it to a natural
isomorphism  
$$
\beta_M: \Hom^{H \ot_R^l H}(H \ot_R^l H,M\ot_R H) \to \Hom_R(H,M),\;
\Phi\mapsto (\id_M\ot_R \counit)\circ \Phi\circ \cop.
$$
An easy computation shows that $\beta_M$ is a morphism of 
$[H\ot_R^l H,-]$-modules, what completes the proof.  
\end{proof}

\noindent {\bf Acknowledgements.} The authors are grateful to Bachuki
Mesablishvili  for helpful advice and hints to the literature. TB would like
to thank Dmitriy Rumynin for drawing his attention to contramodules.
GB expresses her thanks for a 
Bolyai J\'anos Scholarship and support by the Hungarian Scientific Research
Fund OTKA F67910. The work on this paper was partly carried out, when the
first two authors visited the Mathematical Institute  of Polish Academy of
Sciences (IMPAN) in Warsaw, and the support 
of  the European Commission grant MKTD-CT-2004-509794 is acknowledged. 
They would like to thank Piotr M. Hajac for the very warm hospitality in
Warsaw.

\end{document}